\documentclass[12pt,a4paper]{article}
\usepackage{latexsym,amsfonts,amstext,amsbsy,amsmath,amssymb,amscd,graphicx,epsfig,eucal,bbm}
\usepackage{mathrsfs}
\usepackage{pb-diagram}
\usepackage[all]{xy}
\usepackage{makeidx}
\usepackage{indentfirst}
\usepackage{multicol}
\usepackage{graphicx}
\newtheorem{te}{Theorem}[section]

\newtheorem{pr}[te]{Proposition}

\DeclareMathAlphabet{\mathpzc}{OT1}{pzc}{m}{it}
\DeclareMathAlphabet{\filip}{OT1}{rsfs}{m}{it}

\date{}

\begin{document}

\title{The analysis of stochastic stability of stochastic models that describe
 tumor-immune systems}

 \maketitle

\begin{center}
\author{O. Chi\c{s}$^{1}$, A. Sandru$^{2}$, D. Opri\c{s}$^{3}$}\\
\end{center}

\begin{center}
 $ ^{1}$ Euro University "Dr\v{a}gan", Lugoj, Romania,\\ email: chisoana@yahoo.com\\
 $ ^{2}$ Faculty of Exact Sciences, Aurel Vlaicu University, Arad, Romania, \\ email:
  sandruandrea@gmail.com\\
 $ ^{3}$ Faculty of Mathematics and Informatics, West University of Timi\c{s}oara,
  Romania,\\ email: opris@math.uvt.ro
\end{center}

\maketitle \textbf{Abstract}:  In this paper we investigate some
stochastic models for tumor-immune systems. To describe these
models, we used a Wiener process, as the noise has a stabilization
effect. Their dynamics are studied in terms of stochastic
stability in the equilibrium points, by constructing the Lyapunov
exponent, depending on the parameters that describe the model.
Stochastic stability was also proved by constructing a Lyapunov
function. We have studied and and analyzed a Kuznetsov-Taylor like
stochastic model and a Bell stochastic model for tumor-immune
systems. These stochastic models are studied from stability point
of view and they were represented using the second Euler scheme
and Maple 12 software.

\section{Introduction}

Stochastic modeling plays an important role in many branches of
science. In many practical situations perturbations appear and
these are expressed using white noise, modeled by Brownian
motion. We will study stochastic dynamical systems that are used
in medicine, in describing a tumor behavior, but still we don't
know much about the mechanism of destruction and establishment of
a cancer tumor, because a patient may experience tumor regression
and  later a relapse can occur. The need to address not only
preventative measures, but also more successful treatment
strategies is clear. Efforts  along these lines are now being
investigated through immunotherapy (\cite{Ono}, \cite{Vla},
\cite{Whe}).

This tumor-immune study, from theoretical point of view, has been
done for two cell populations:  effector cells and tumor cells. It
was predicted a threshold above which there is uncontrollable
tumor growth, and below which the disease is attenuated with
periodic exacerbations occurring every 3-4 months. There was also
shown that the model does have stable spirals, but the
Dulac-Bendixson criterion demonstrates that there are no stable
closed orbits. It is consider ODE's for the populations of immune
and tumor cells and it is shown that survival increases if the
immune system is stimulated, but in some cases an increase in
effector cells increases the chance of tumor survival.

In the last years, stochastic growth models for cancer cells were
developed, we mention the papers of W.Y. Tan and C.W. Chen
\cite{Tan}, N. Komarova, G. Albano and V.Giorno \cite{Albano}, L.
Ferrante, S. Bompadre, L. Possati and L. Leone \cite{Ferr}, A.
Boondirek Y. Lenbury, J. Wong-Ekkabut, W. Triampo, I.M. Tang, P.
Picha \cite{Boo}.

Our goal in this paper is to  construct stochastic models and to
analyze their behavior around the equilibrium point. In these
points stability is studied by analyzing the Lyapunov exponent,
depending of the parameters of the models. Numerical simulations
are done using a deterministic algorithm with an ergodic invariant
measure. In this paper the authors studied and analyzed two
stochastic models. In Section 2, we considered a Kuznetsov and
Taylor stochastic model. Beginning from the classical one, we have
studied the case of positive immune response. We gave the
stochastic model and we analyzed it in the equilibrium points.
Numerical simulations for this new model are presented in Section
2.1. In Section 3 we presented a general family of tumor-immune
stochastic systems and from this general representation we
analyzed Bell model. We wrote this model as a stochastic model,
using Annexe 1, and we discussed its behavior around the
equilibrium points. We have proved stochastic stability around
equilibrium point using two methods. The first one consists of
expressing the Lyapunov exponent, and then drawing the conclusion when the
considered system is stable. The second method is a way of
constructing a Lyapunov function and determining sufficient
conditions such that the system is stable. Numerical simulations
were done using the software Maple 12 and we implemented the
second order Euler scheme for a representation of the discussed
stochastic models.

\section{Kuznetsov and Taylor stochastic model}

The study of tumor-immune interaction is determined by the behavior of two
po\-pulations of cells: effector cells and tumor cells. We will
construct the stochastic models using well known deterministic
models and we analyze stochastic stability around the equilibrium
points. The analysis is done using Lyapunov exponent method.

 We will begin our study from the deterministic model
of Kuznetsov and Taylor \cite{Kuz}. This model describes the
response of effector cells to the growth of tumor cells and takes
into consideration the penetration of tumor cells by effector
cells, that causes the interaction of effector cells. This model
can be represented in the following way:
\begin{equation}\label{3}
\left \{%
\begin{array}{ll}
\dot{x}(t)=a_1 - a_2 x(t)+ a_3x(t)y(t),\\
\dot{y}(t)= b_1y(t)(1 - b_2y(t)) - x(t)y(t),\\
\end{array}%
\right.
\end{equation}
where initial conditions are $x(0)=x_0 >0,\, y(0)=y_0 >0$ and
$a_3$ is the immune response to the appearance of the tumor cells.

In this paper we consider the case of $a_3>0,$ that means that
immune response is positive. For the equilibrium states $P_1$ and
$P_2,$ we study the asymptotic behavior  with respect to the
parameter $a_1$ in (\ref{3}). For $b_1a_2<a_1,$ the system
(\ref{3}) has the equilibrium states $P_1(x_1,y_1)$ and
$P_2(x_2,y_2),$ with
\begin{equation}\label{4}
x_1=\frac{a_1}{a_2}, \, y_1=0,
\end{equation}
\begin{equation}\label{5}
x_2=(b_1(a_3-b_2a_2)+\sqrt(\Delta))/(2a_3), \,
y_2=(b_1(a_3+b_2a_2)-\sqrt(\Delta)/(2b_1b_2a_3)
\end{equation}
where $\Delta =b_1^2(b_2a_2-a_3)^2+4b_1b_2a_1a_3.$\\

We associate a stochastic system of differential equations to
the ordinary system of differential equations (\ref{3}).

In \cite{Kuz} it is shown that there is an $a_{10},$ such that if
$0<a_1<a_{10},$ then the equilibrium state $P_1$ is asymptotically
stable, and for $a_1>a_{10}$ the equilibrium state $P_1$ is
unstable. If $a_1<a_{10},$ then the equilibrium state $P_2$ is
unstable and for $a_1>a_{10}$ it is asymptotically stable.

\par Let us consider $(\Omega,\mathcal{F}_{t\geq 0},\mathcal{P} )$
a filtered probability space and $(W(t))_{t\geq 0}$ a standard
Wiener process adapted to the filtration $(\mathcal{F_t})_{t\geq
0}.$ Let $\{X(t,\omega)=(x(t),y(t))\}_{t\geq 0}$ be a stochastic
process.

The system of It\^o equations associated to system (\ref{3}) is
given by
\begin{equation}\label{6}
\left \{%
\begin{array}{ll}
x(t)=x_0+\int_0^t (a_1-a_2x(s)+a_3x(s)y(s))ds+\int_0^t
g_1(x(s),y(s))dW(s),\\
 y(t)=y_0+\int_0^t
((b_1y(s)(1-b_2y(s))-x(s)y(s))ds+\int_0^t
g_2(x(s),y(s))dW(s),\\
\end{array}
\right.
\end{equation}where the first integral is a Riemann  integral, and
the second one is an It\^o integral. $\{W(t)\}_{t>0}$ is a Wiener
process \cite{Schu}.

The functions $g_1(x(t),y(t))$ and $g_2(x(t),y(t))$ are given in
the case when we are working in the equilibrium state. In $P_1$
those functions have the following form
\begin{equation}\label{7}
\begin{array}{ll}
g_1(x(t),y(t))=b_{11}x(t)+b_{12}y(t)+c_{11},\\
\quad \\
g_2(x(t),y(t))=b_{21}x(t)+b_{22}y(t)+c_{21},\\
\end{array}
\end{equation}where
\begin{equation}\label{8}
c_{11}=-b_{11}x_1-b_{12}y_1, \, c_{21}=-b_{21}x_1-b_{22}y_1.
\end{equation}

In the equilibrium state $P_2,$ the functions $g_1(x(t),y(t))$ and
$g_2(x(t),y(t))$ are given by
\begin{equation}\label{9}
\begin{array}{ll}
g_1(x(t),y(t))=b_{11}x(t)+b_{12}y(t)+c_{12},\\
\quad \\
g_2(x(t),y(t))=b_{21}x(t)+b_{22}y(t)+c_{22},\\
\end{array}
\end{equation}where
\begin{equation}\label{8}
c_{12}=-b_{11}x_2-b_{12}y_2, \, c_{22}=-b_{21}x_2-b_{22}y_2.
\end{equation}

The functions $g_1(x(t),y(t))$ and $g_2(x(t),y(t))$ represent the
volatilisations of the stochastic equations and they are the
therapy test functions.

\subsection{The analysis of SDE (\ref{6}). Numerical simulation.}

Using the formulae from Annexe 1, Annexe 2, and  Maple 12
software, we get the following results, illustrated in the below
figures. For numerical simulations we use the following values of
parameters:
$$a_1=0.1181, \, a_2=0.3747, \, a_3=0.01184, \, b_1=1.636, \, b_2=0.002.$$
The matrices $A$ and $B$ are given, in the equilibrium point
$P_1(\frac{a_1}{a_2},0)$ by
$$A=\begin{pmatrix}
-a_2+a_3y_1 & a_3x_1 \\
-y_1 & b_1-2b_2y_1-x_1\\
\end{pmatrix}, \quad
B=\begin{pmatrix}
10 & -2 \\
2 & 10\\
\end{pmatrix}.$$
In a similar way, matrices $A$ and $B$ are defined in the other
equilibrium point
$$P_2\Big(\frac{(-b_1(b_2a_2-a_3)+\sqrt{\Delta})}{2a_3},
\frac{(b_1(b_2a_2+a_3)-\sqrt{\Delta})}{2b_1b_2a_3}\Big),$$ with
$\Delta=b_1^2(b_2a_2-a_3)^2+4b_1b_2a_1a_3.$

Using second order Euler scheme, for the ODE system (\ref{3}) and SDE system (\ref{6}), we get the following orbits presented in the figures above.

\begin{center}\begin{tabular}{cc}
\epsfxsize=6cm \epsfysize=5cm
 \epsffile{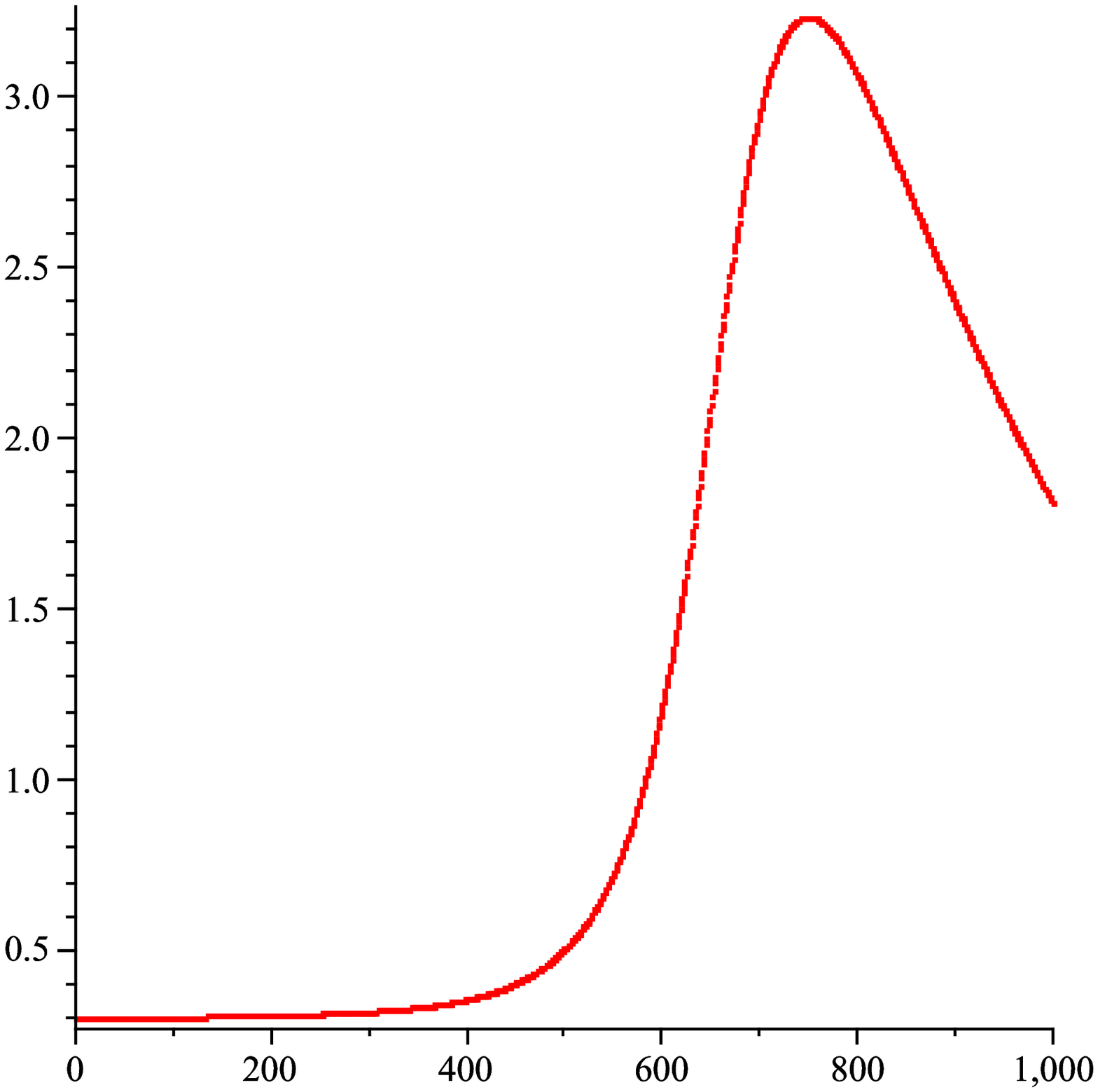}  &
\epsfxsize=6cm \epsfysize=5cm
\epsffile{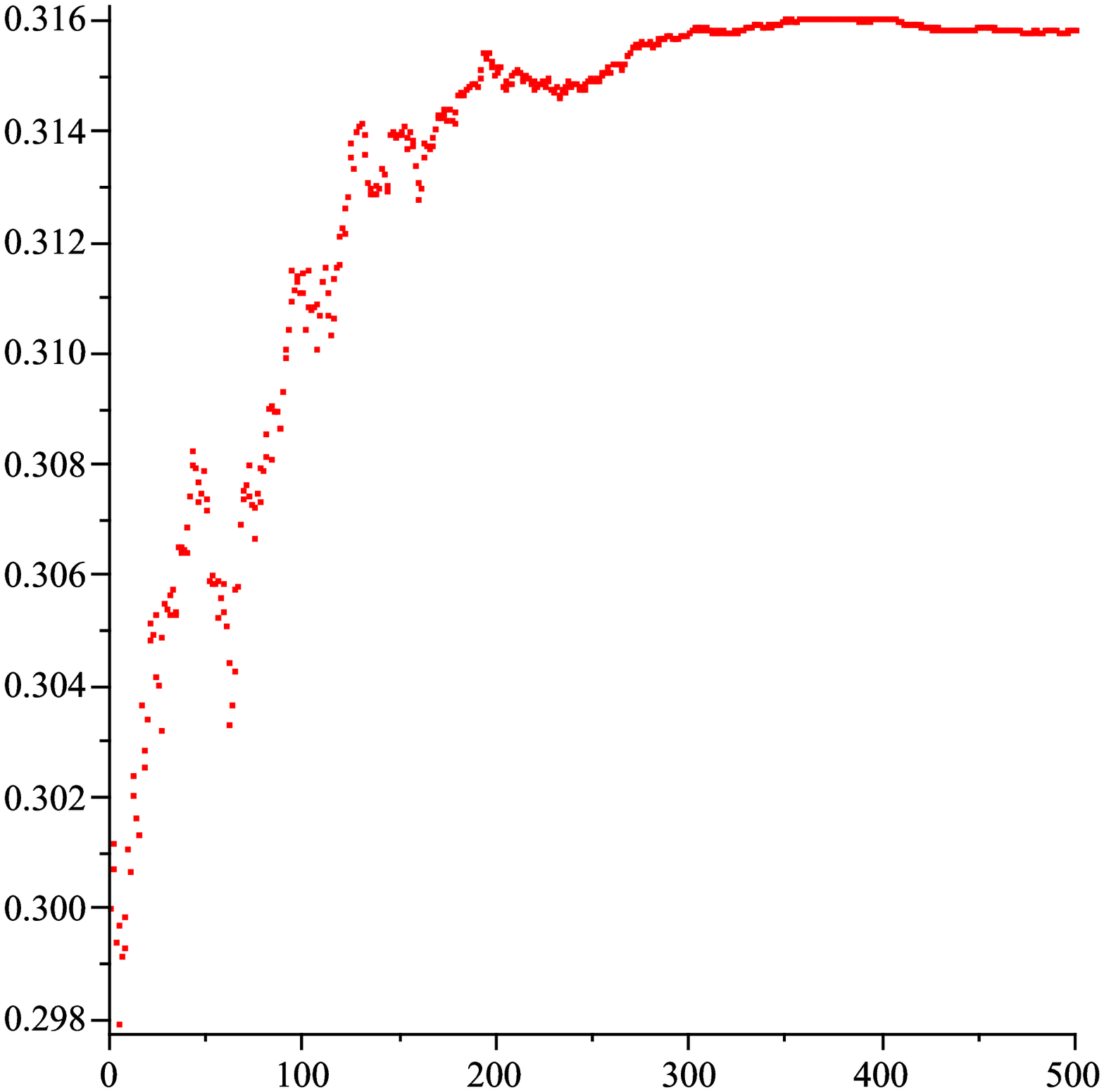} \\
  \footnotesize Figure 1: $(n,x(n))$ in $P_1$ for ODE (\ref{3})  &  \footnotesize Figure 2: $(n,x(n,\omega))$ in $P_1$  for SDE (\ref{6})\\
   &   \\
        &   \\
\epsfxsize=6cm \epsfysize=5cm
 \epsffile{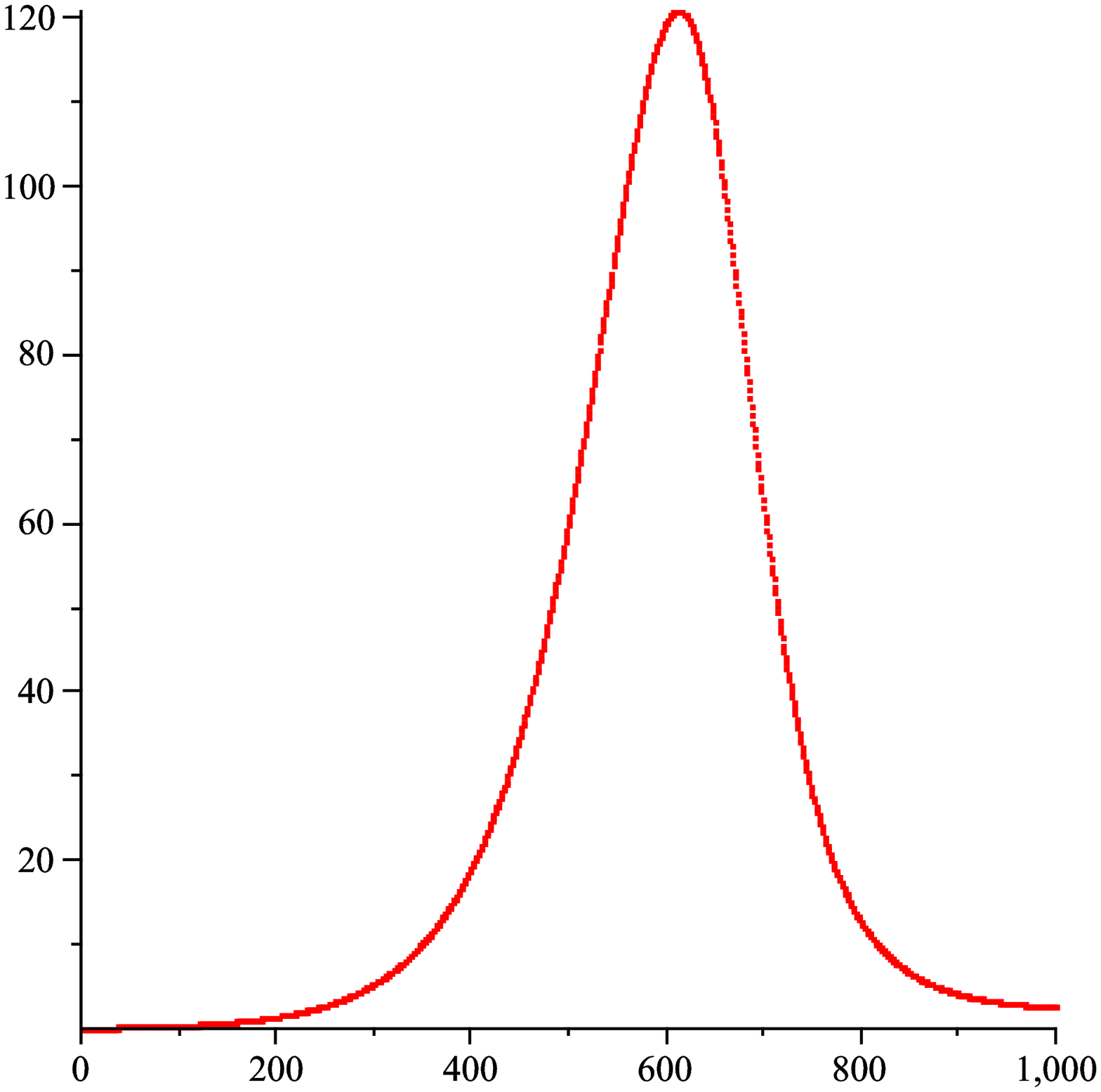}  &
\epsfxsize=6cm \epsfysize=5cm
\epsffile{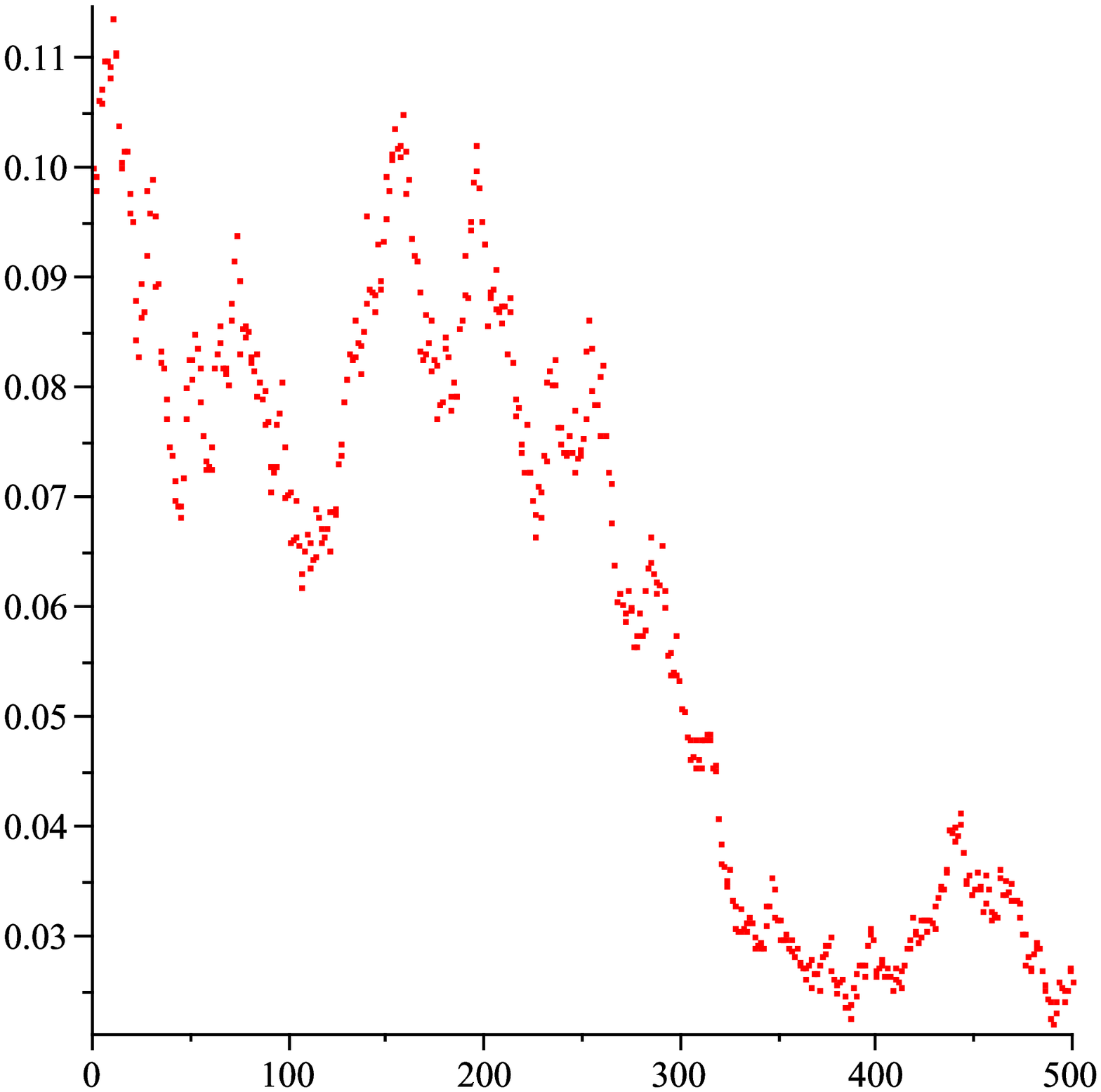} \\
    \footnotesize Figure 3: $(n,y(n))$ in $P_1$  for ODE (\ref{3}) &  \footnotesize Figure 4: $(n,y(n,\omega))$ in $P_1$  for SDE (\ref{6})\\
    &   \\
        &   \\

        \end{tabular}
\end{center}
\begin{center}

\begin{tabular}{cc}
\epsfxsize=6cm \epsfysize=5cm
 \epsffile{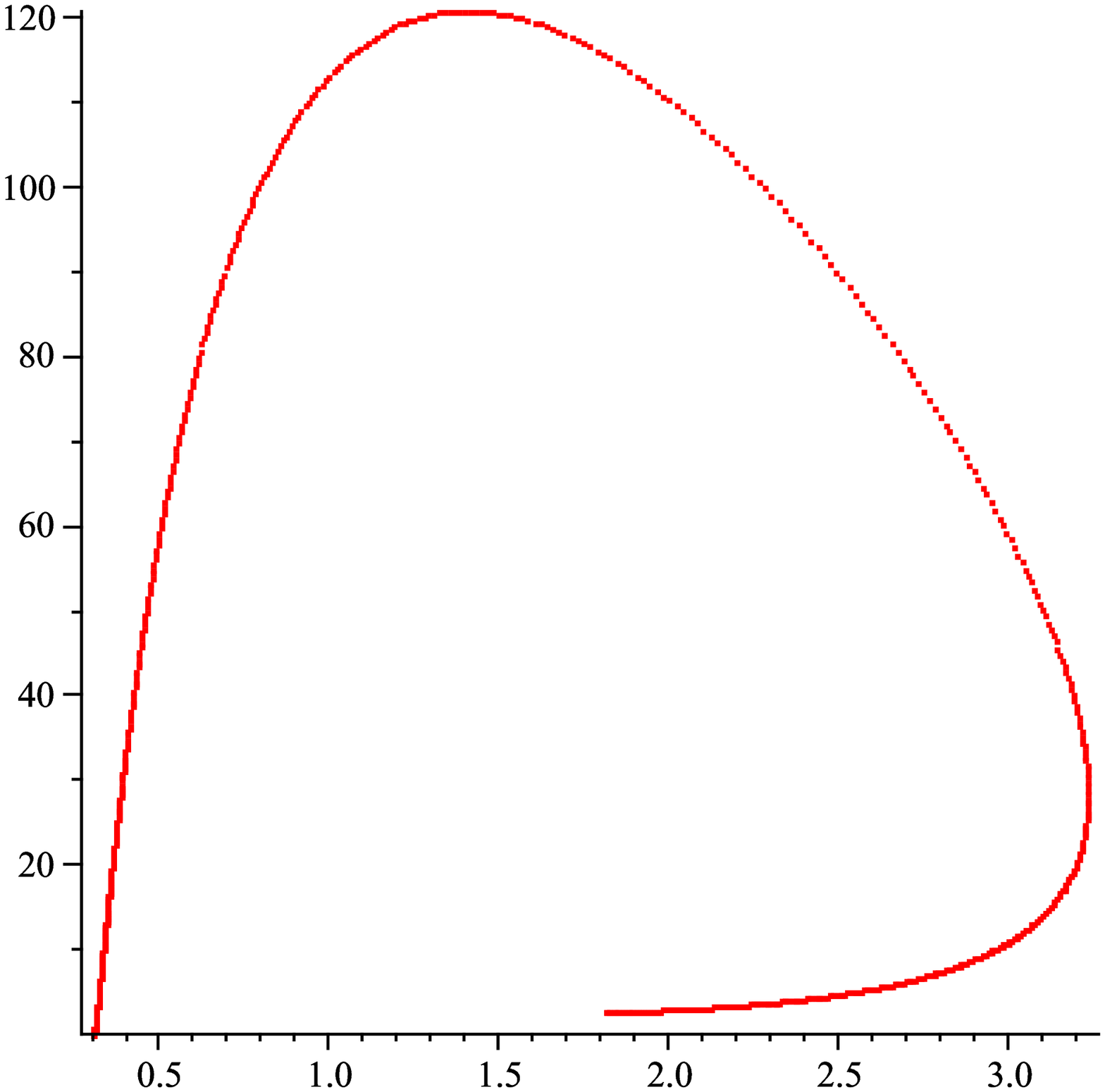}  &
\epsfxsize=6cm \epsfysize=5cm
\epsffile{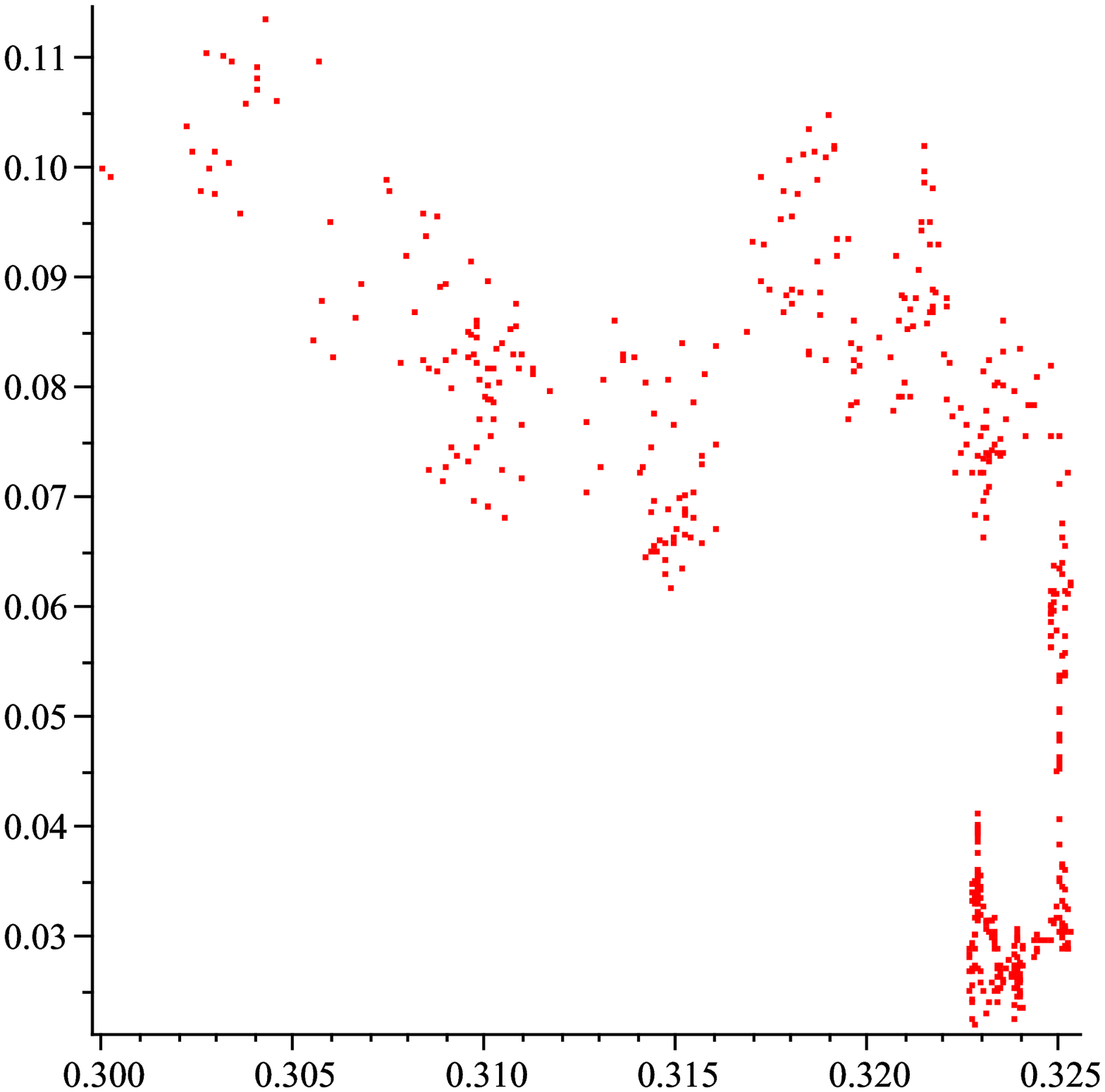} \\
  \footnotesize  Figure 5: $(x(n),y(n))$  in $P_1$  for ODE (\ref{3})&  \footnotesize Figure 6: $(x(n,\omega),y(n,\omega))$ in $P_1$  for SDE (\ref{6})\\
      &   \\
        &   \\ \end{tabular}
\end{center}
\begin{center}

\begin{tabular}{cc}
\epsfxsize=6cm \epsfysize=5cm
 \epsffile{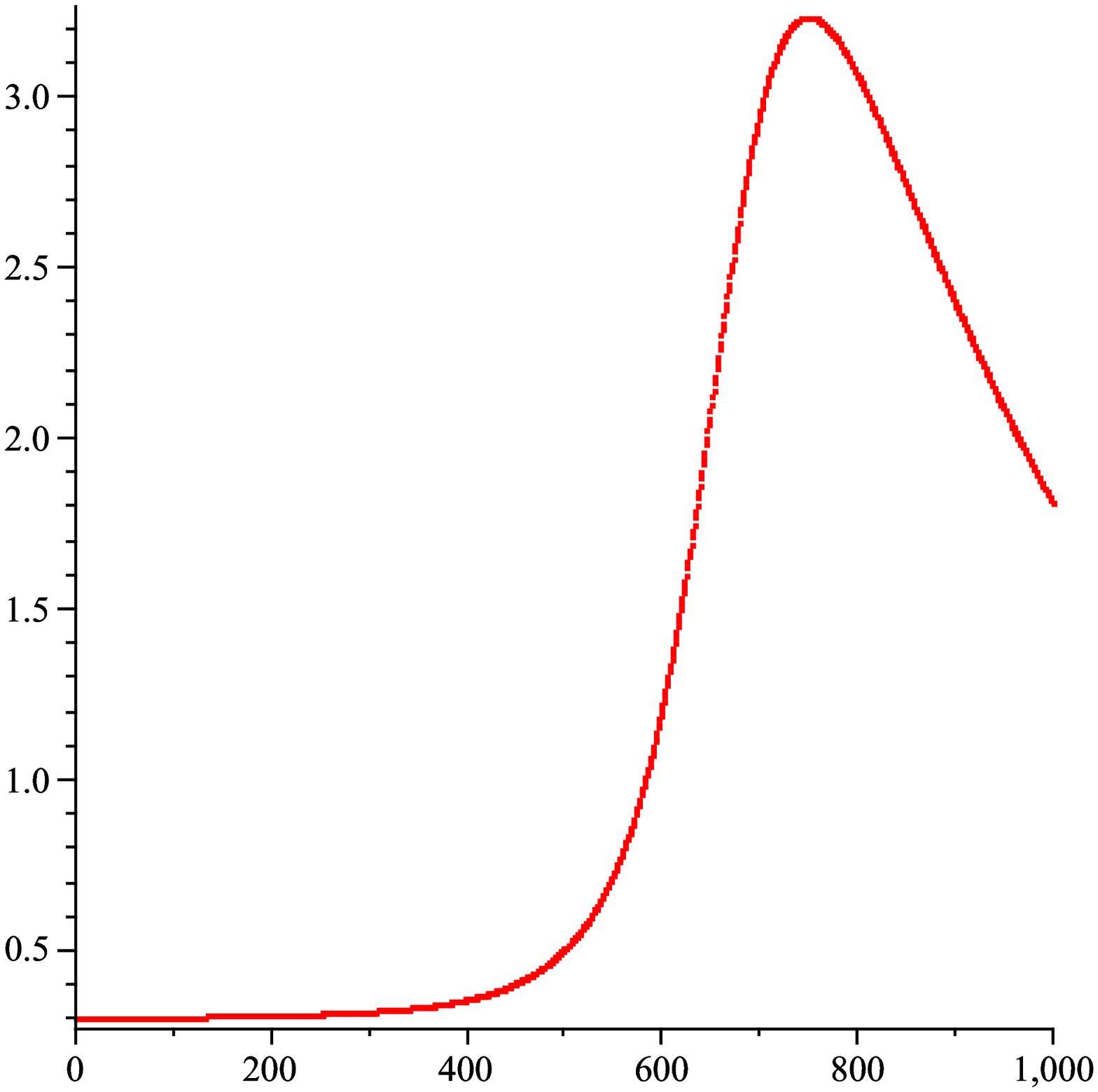}  &
\epsfxsize=6cm \epsfysize=5cm
\epsffile{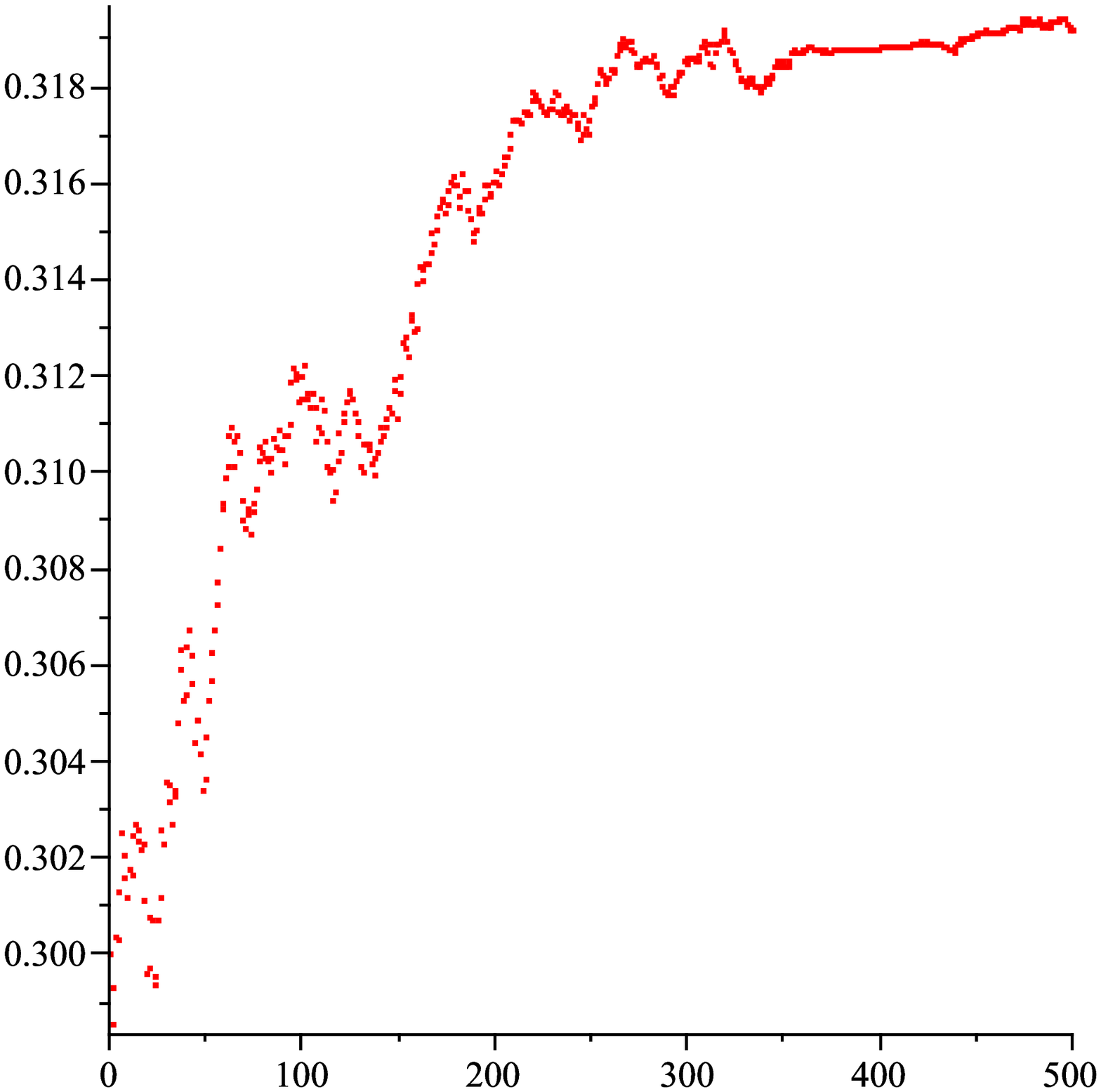} \\
    \footnotesize Figure 7: $(n,x(n))$  in $P_2$  for ODE (\ref{3})&  \footnotesize Figure 8: $(n,x(n,\omega))$ in $P_2$  for SDE (\ref{6})\\
       &   \\
        &   \\

 \end{tabular}
\end{center}
\begin{center}

\begin{tabular}{cc}

        \epsfxsize=6cm \epsfysize=5cm
 \epsffile{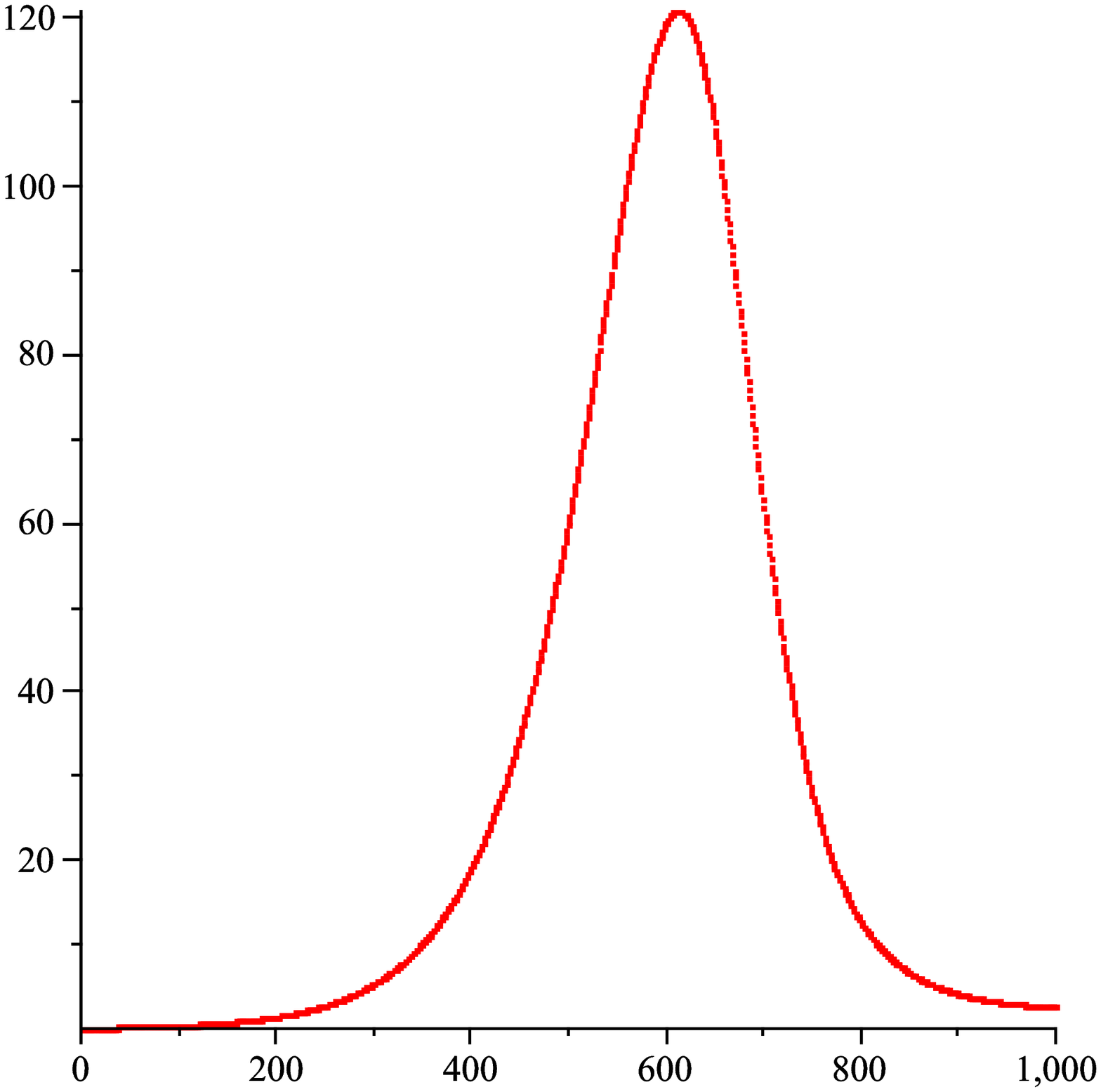}  &
\epsfxsize=6cm \epsfysize=5cm
\epsffile{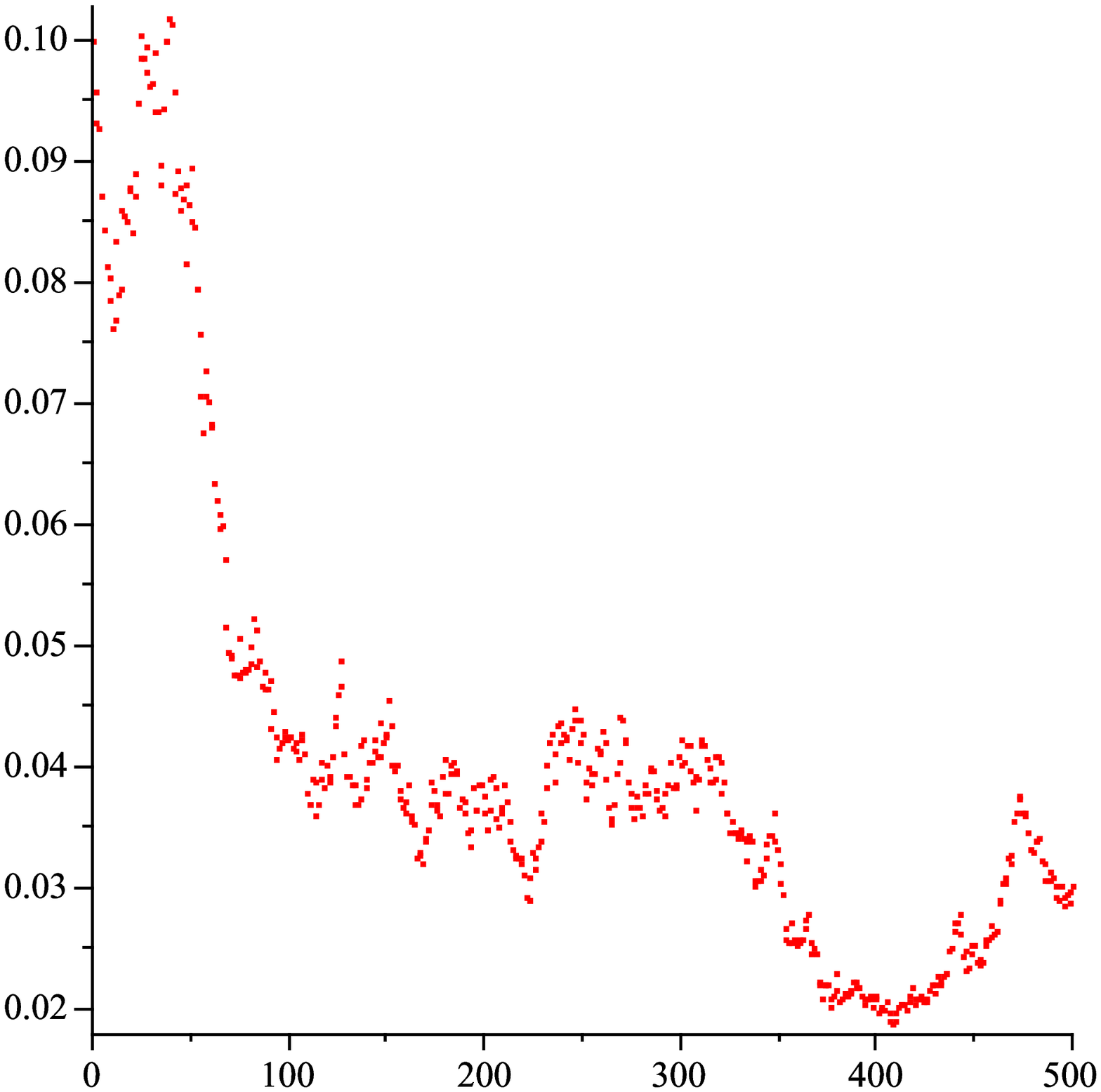} \\
   \footnotesize Figure 9: $(n,y(n))$ in $P_2$   for ODE (\ref{3})&  \footnotesize Figure 10: $(n,y(n,\omega))$ in $P_2$  for SDE (\ref{6})\\
        &   \\
        &   \\
         \end{tabular}
\end{center}
\begin{center}

\begin{tabular}{cc}

\epsfxsize=6cm \epsfysize=5cm
 \epsffile{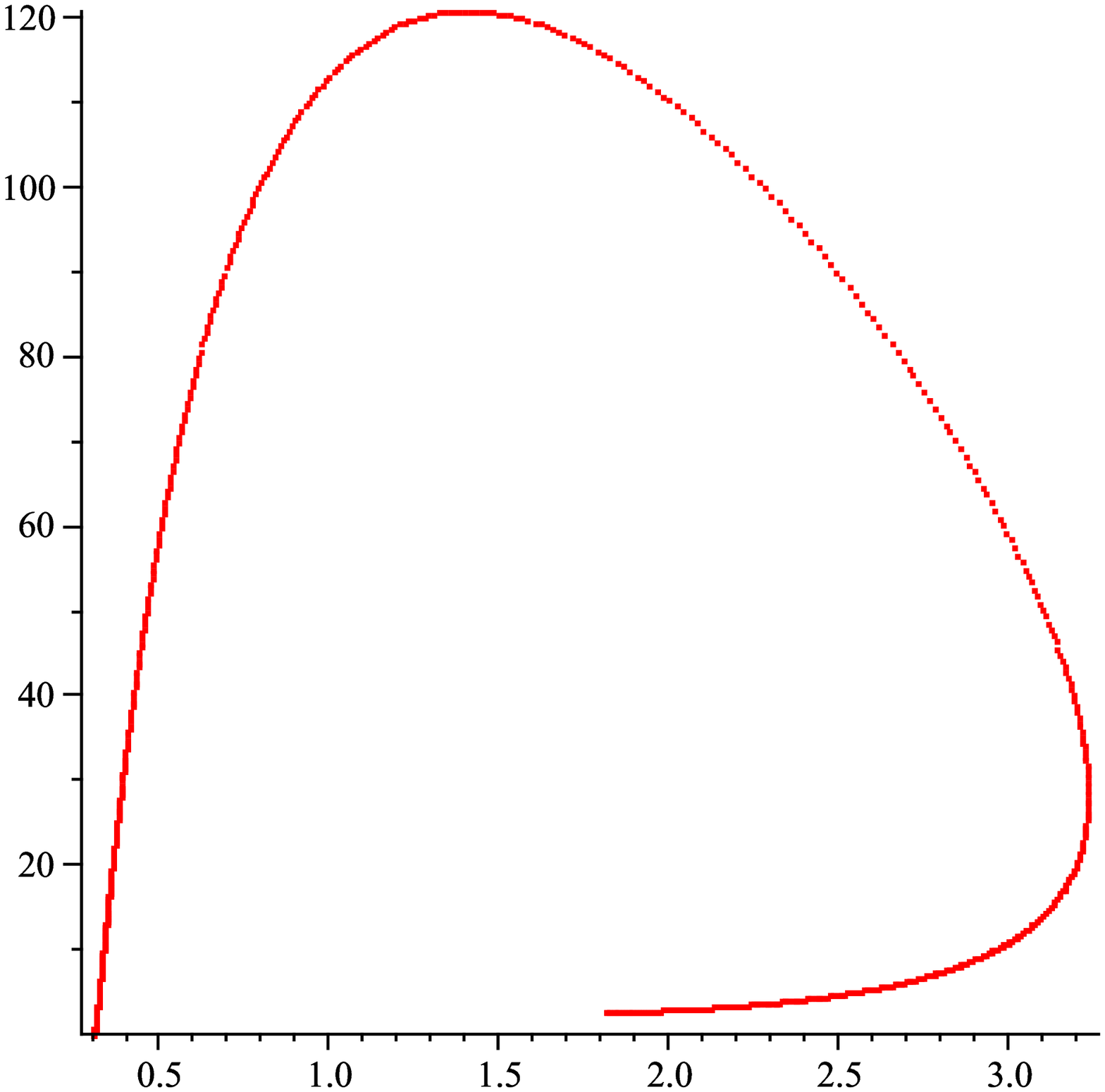}  &
\epsfxsize=6cm \epsfysize=5cm
\epsffile{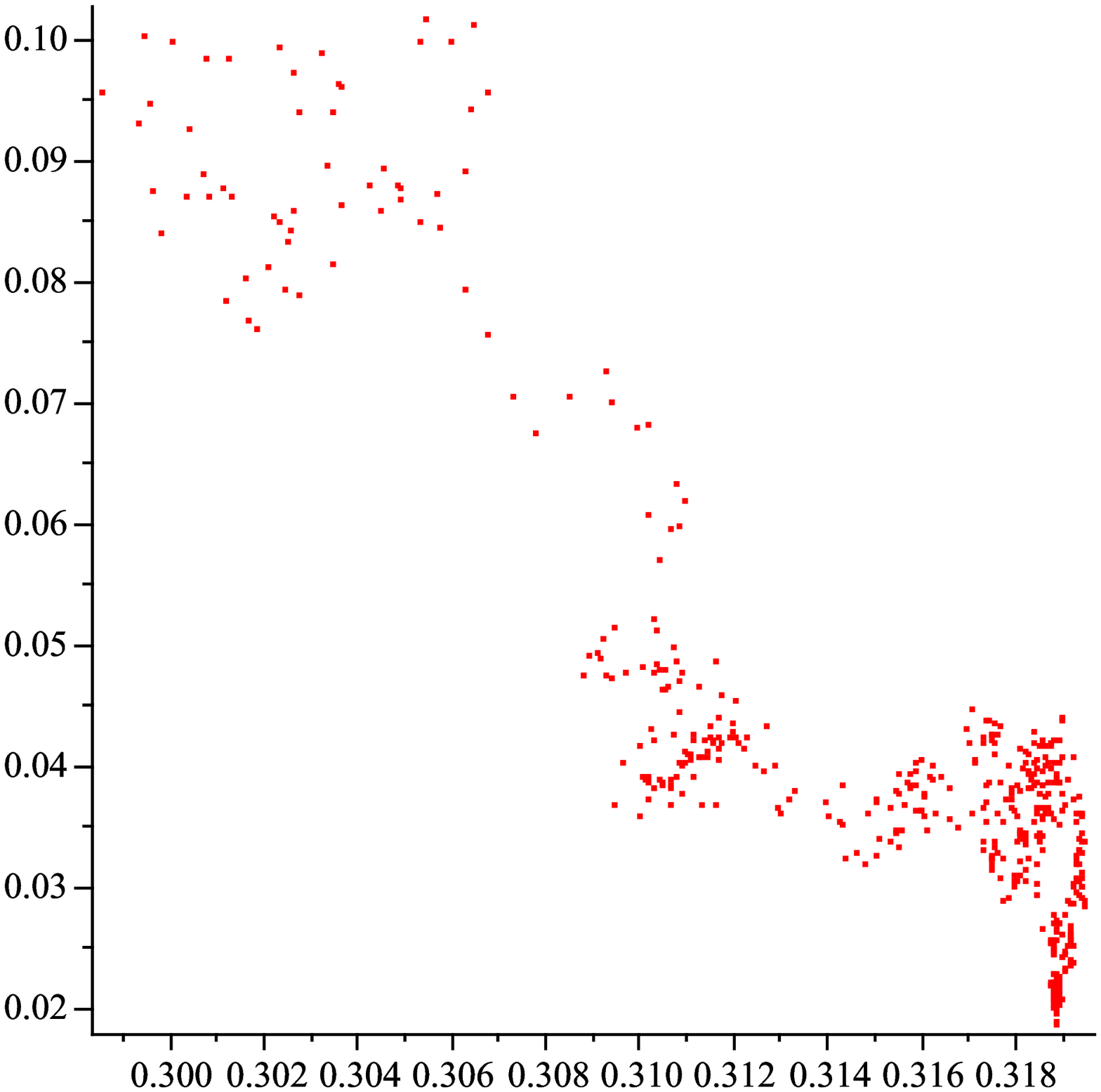} \\
 \footnotesize Figure 11: $(x(n),y(n))$ in $P_2$  for ODE (\ref{3})&  \footnotesize Figure 12: $(x(n,\omega),y(n,\omega))$ in $P_2$  for SDE (\ref{6})\\
          &   \\
        &   \\
         \end{tabular}
\end{center}
\begin{center}

The Lyapunov exponent variation, with the variable parameter $b_{11}=\alpha,$ is given in Figure 13 for $P_1,$ and in Figure 14 for the equilibrium point $P_2.$

\begin{tabular}{cc}

\epsfxsize=6cm \epsfysize=5cm
 \epsffile{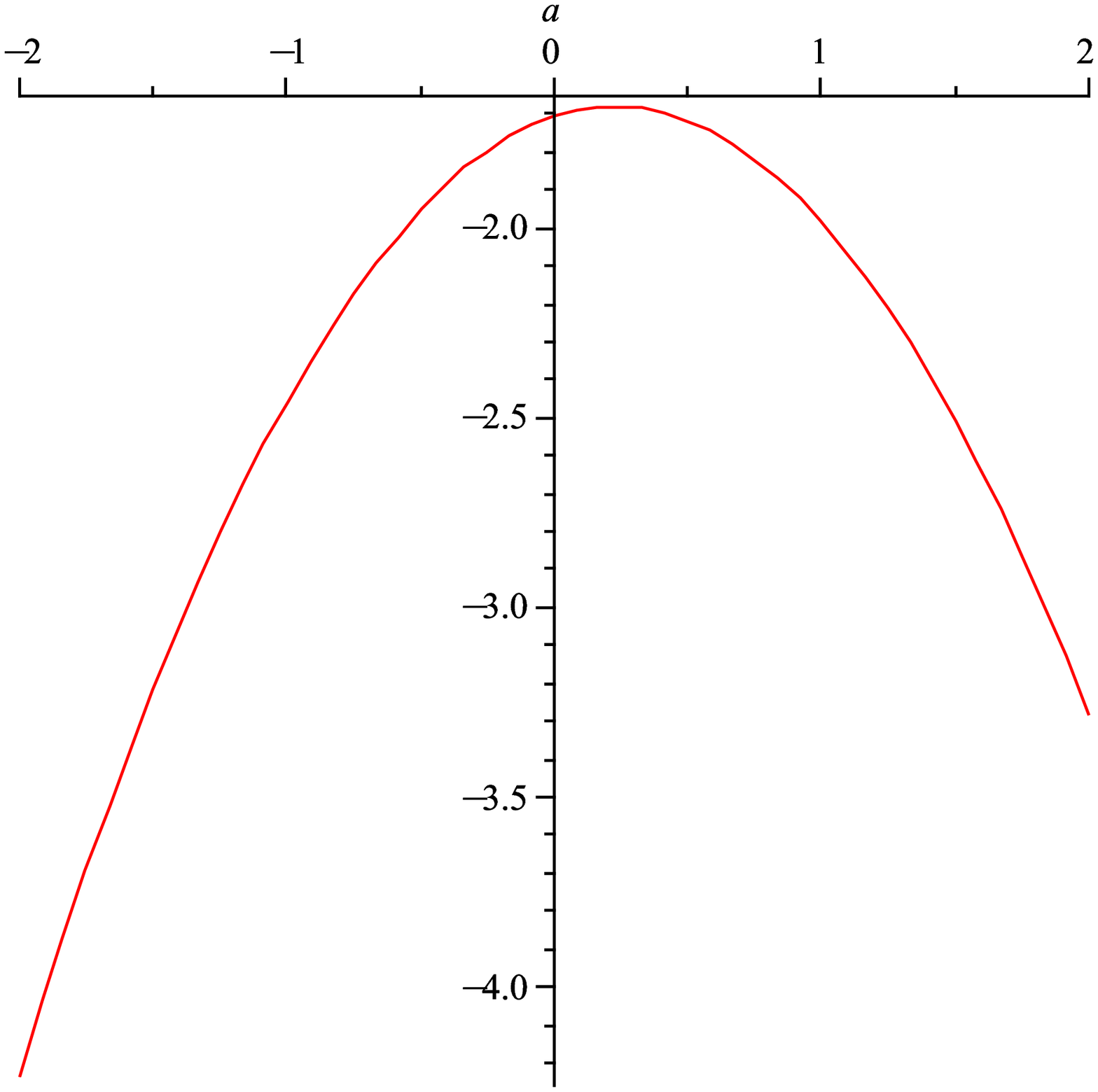}  &
\epsfxsize=6cm \epsfysize=5cm
\epsffile{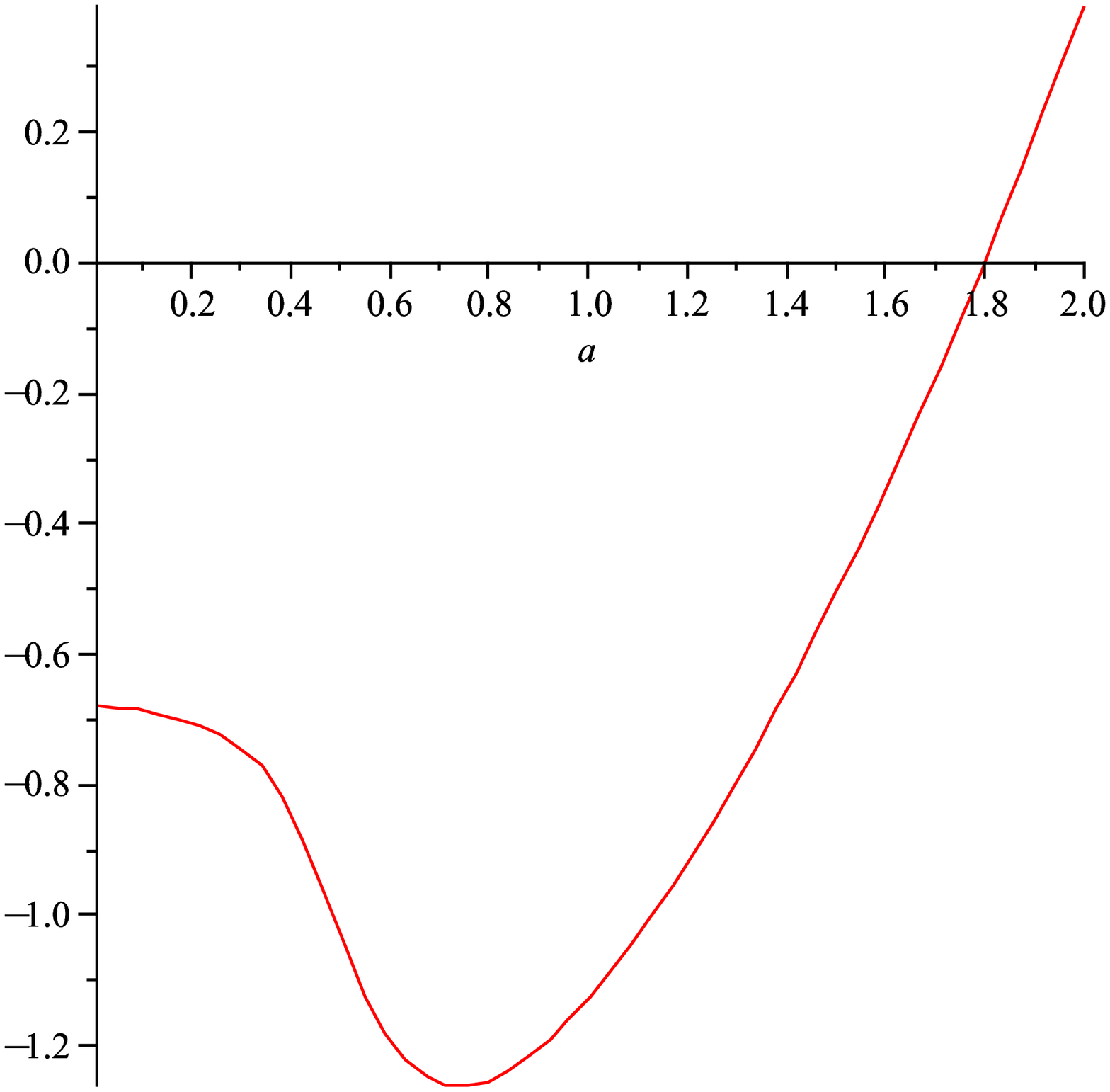} \\
  Figure 13: $(\alpha,\lambda(\alpha))$ in $P_1$  for ODE (\ref{3})  &  Figure 14: $(\alpha,\lambda(\alpha))$ in $P_2$  for SDE (\ref{6})\\
        &   \\
        &   \\
\end{tabular}
\end{center}

The Lyapunov exponent, for the equilibrium point $P_1$ is
negative, so $P_1$ is asymptotically stable for each
$\alpha\in\mathbb{R}.$ For the second equilibrium point $P_2,$
this point is asymptotically stable for all values when
$\lambda<0,$ that means that $P_2$ is unstable for
all $\alpha\in (-\infty,-1.8)\cup (1.8,\infty). $

\section{A general family of tumor-immune stochastic systems}

A Volterra-like model was proposed in \cite{Soto} for the
interaction between a population of tumor cells (whose number is
denoted by $x$) and a population of lymphocyte cells ($y$) and it
is given by
\begin{equation}\label{g1}
\left \{%
\begin{array}{ll}
\dot{x}(t)=ax(t)-bx(t)y(t),\\
\dot{y}(t)=dx(t)y(t)-fx(t)-kx(t)+u,\\
\end{array}%
\right.
\end{equation}
where the tumor cells are supposed to be in exponential growth
(which is, however, a good approximation only for the initial
phases of the growth) and the presence of tumor cells implies a
decrease of the "input rate" of lymphocytes.

A general representation for such models can be considered in the
form given by d'Onofrio in \cite{Ono}:
\begin{equation}\label{g2}
\left \{%
\begin{array}{ll}
\dot{x}(t)=f_1(x(t),y(t)), \, \dot{y}(t)=f_2(x(t),y(t)),\\
x(0)=x_0, \, y(0)=y_0,
\end{array}%
\right.
\end{equation}where $x$ is the number of tumor cells, $y$ the number of effector cells
     of immune system and
\begin{equation}
\begin{array}{ll}
     f_1(x,y)=x(h_1x-h_2xy),\\
     \quad \\
     f_2(x,y)=(h_3x-h_4x)y+h_5x.
     \end{array}
     \end{equation}
The functions $h_1,h_2,h_3,h_4,h_5$ are given such that the system
(\ref{g2}) admits the equilibrium point $P_1(x_1,y_1),$ with
$x_1=0, \, y_1>0,$ and $P_2(x_2,y_2),$ with $x_1\neq0, \, y_2>0.$

Deterministic models of this general form are the following
\begin{description}
\item[Volterra model \cite{Volt}] if $h_1(x)=a_1, \, h_2(x)=a_2x,
\, h_3(x)=b_3x, \, h_4(x)=b_2 $ and $h_5(x)=-b_1x;$
    \item[Bell model]  $h_1(x)=a_1x, \, h_2(x)=a_2x, \,
    h_3(x)=b_1x,$
$h_4(x)=b_3 $ and $h_5(x)=-b_2x+b_4;$
    \item[Stepanova model \cite{Ste}] with
    $h_1(x)=a_1, \, h_2(x)=1, \, h_3(x)=b_1 x, \, h_4(x)=b$ and
$h_5(x)=-b_2x+b_4;$
    \item[Vladar-Gonzalez model \cite{Vla}] if in (\ref{g2}) we consider
   $h_1(x)= \log(K/x),$ $ h_2(x)=1, \, h_3(x)=b_1 x, \, h_4(x)=
   b_2+b_3 x^2$ and
$h_5(x)=1;$
    \item[Exponential model \cite{Whe}] if in (\ref{g2}) we
    consider $h_1(x)=1,$ $h_2(x)=1,$ $h_3(x)=b_1x,$
    $h_4(x)=b_2+b_3x^2,$ and $h_5(x)=1;$
    \item[Logistic model \cite{Maru}] if in (\ref{g2}) we consider
   $h_1(x)=1-\frac{a_1}{x},$ $h_2(x)=1,$ $h_3(x)=b_1x,$
    $h_4(x)=b_2+b_3x^2,$ and $h_5(x)=1.$
\end{description}
The analysis of these models was proven also using numerical
simulations.

For a considered filtered probability space
$(\Omega,\mathcal{F}_{t\geq0},\mathcal{P} )$ and a standard Wiener
process $(W(t))_{t\geq0},$ we consider the stochastic process in two
dimensional space $(\mathcal{F_t})_{t\geq0}.$

The system of It\^o equations associated to system (\ref{g2}) is
given, in the equilibrium point $P(x_0,y_0),$ by
\begin{equation}\label{g3}
\left \{%
\begin{array}{ll}
x(t)=x_0+\int_0^t [x(s)(h_1(x(s))-h_2(x(s))y(s)]ds+\int_0^t
g_1(x(s),y(s))dW(s),\\
y(t)=y_0+\int_0^t [(h_3(x(s))-h_4(x(s)))y(s)+h_5(x(s))]ds+\int_0^t
g_2(x(s),\\
\quad  \quad \quad y(s))dW(s),\\
\end{array}
\right.
\end{equation}where the first integral is a Riemann  integral,
the second one is an It\^o integral and $\{W(t)\}_{t>0}$ is a Wiener
process \cite{Schu}.

The functions $g_1(x,y)$ and $g_2(x,y)$ are given in the case when
we are working in the equilibrium state $P_e$, and they are given
by
\begin{equation}
\begin{array}{ll}
g_1(x,y)=b_{11}x+b_{12}y+c_{1e},\\
\quad \\
g_2(x,y)=b_{21}x+b_{22}y+c_{2e},\\
\end{array}
\end{equation}where
\begin{equation}\label{8}
c_{ie}=-b_{i1}x_e-b_{i2}y_e, \, i=1,2,
\end{equation}and $b_{ij}\in \mathbb{R}, \, i,j=1,2.$

\subsection{Analysis of Bell model. Numerical simulations.}

\subsubsection{Lyapunov exponent method}

Following the algorithm for determining the Lyapunov exponent
(Annexe 1) and the description of second order Euler scheme (Annexe 2) in Maple 12 software, we get the following
results, illustrated in the figures below. For numerical
simulations we use the following values of parameters:
$$a_1=2.5, \, a_2=1, \, b_1=1, \, b_2=0.4,\, b_3=0.95, \, b_4=2.$$
The matrices $A$ and $B$ are given, in the equilibrium point
$P_1\Big(0,\frac{b_4}{b_3}\Big)$ by
$$A=\begin{pmatrix}
-a_2y_1+a_1 & -a_2x_1 \\
 -b_2+b_1y_1 & b_1x_1-b_3\\
\end{pmatrix}, \quad
B=\begin{pmatrix}
\alpha & -\beta \\
\beta & \alpha\\
\end{pmatrix},$$with $\alpha=a\in \mathbb{R}, \, \beta=-2.$ In a similar way the
matrices $A$ and $B$ are defined in the other equilibrium point
$P_2(\frac{a_1b_3-a_2b_4}{a_1b_1-a_2b_2},\frac{a_1}{a_2}).$

\begin{center}\begin{tabular}{cc}
\epsfxsize=6cm \epsfysize=5cm
 \epsffile{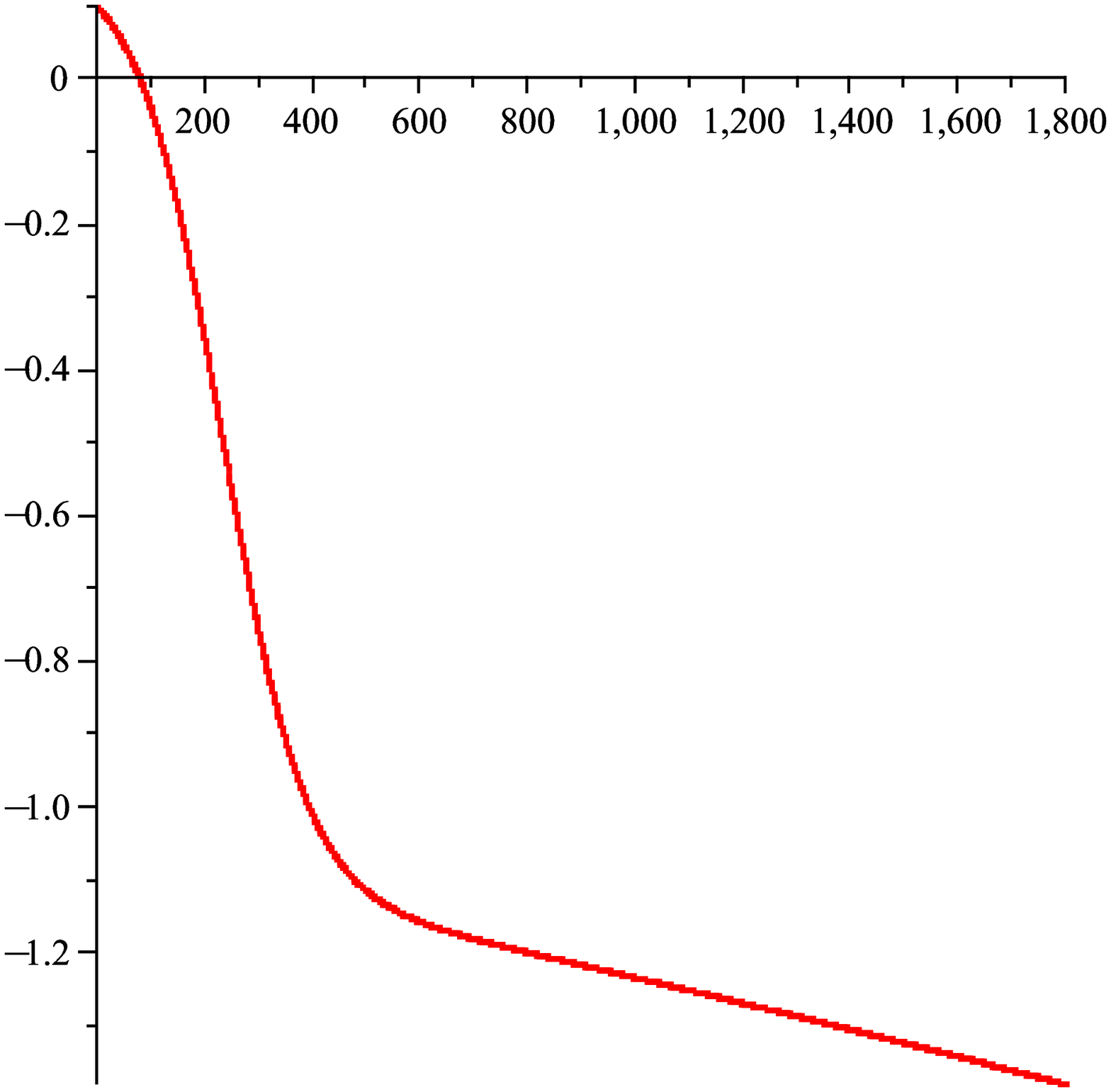}  &
\epsfxsize=6cm \epsfysize=5cm
\epsffile{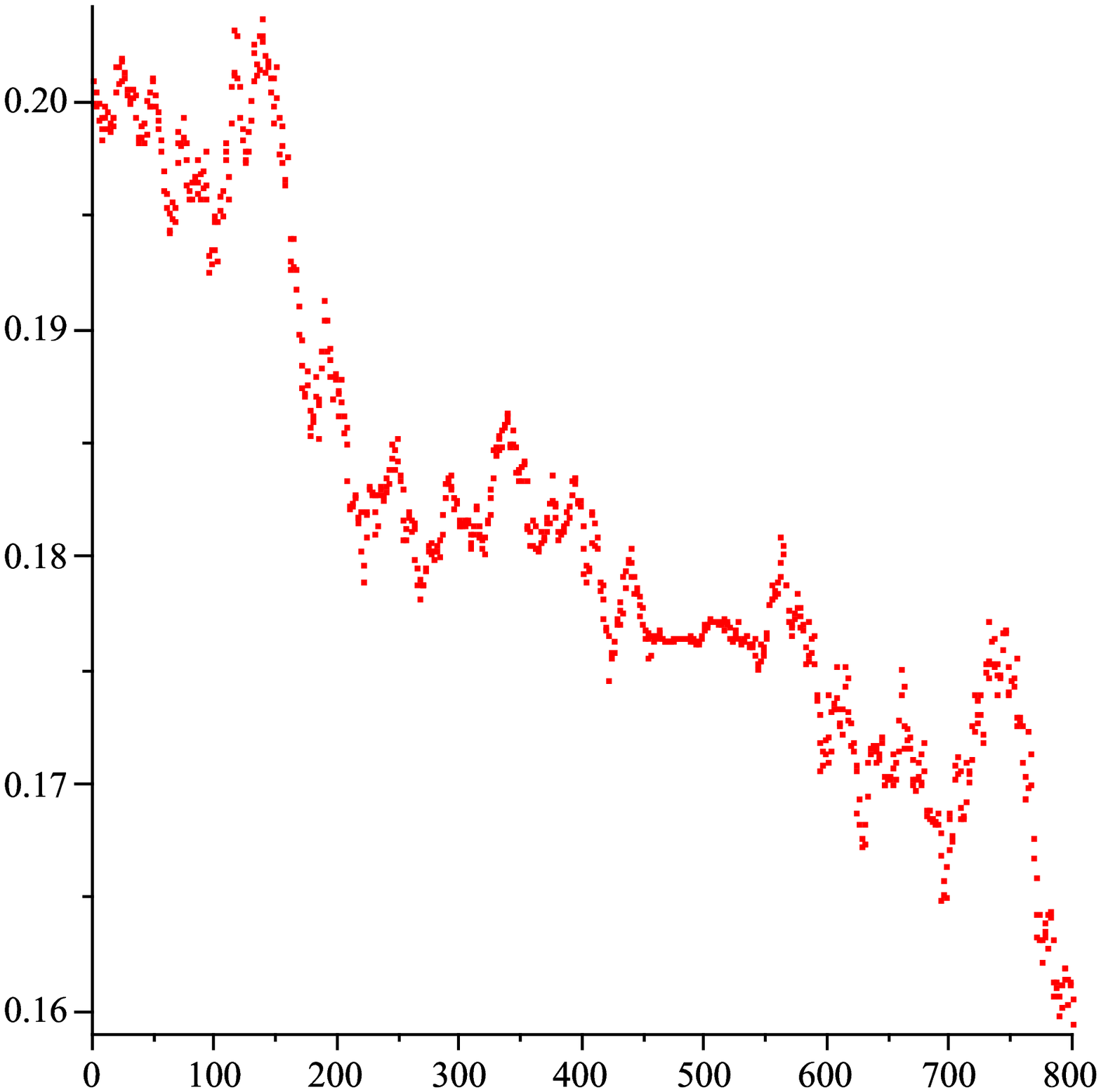} \\
 \footnotesize Figure 15: $(n,x(n))$  in $P_1$ for ODE (\ref{g2}) &\footnotesize Figure 16: $(n,x(n,\omega))$ in $P_1$  for SDE (\ref{g3})\\
 &   \\
        &   \\
        \end{tabular}
\end{center}

\begin{center}
\begin{tabular}{cc}
\epsfxsize=6cm \epsfysize=5cm
 \epsffile{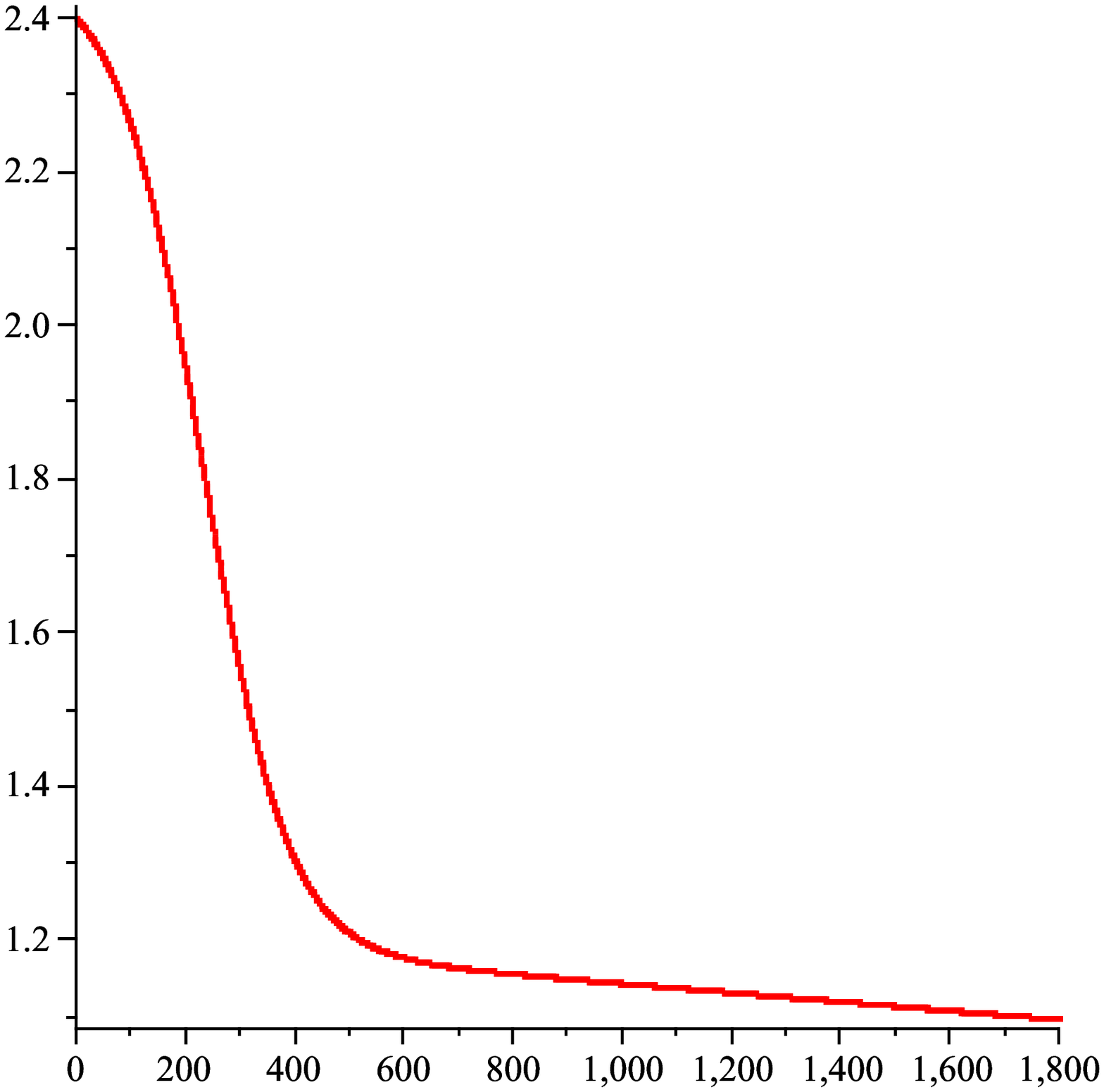}  &
\epsfxsize=6cm \epsfysize=5cm
\epsffile{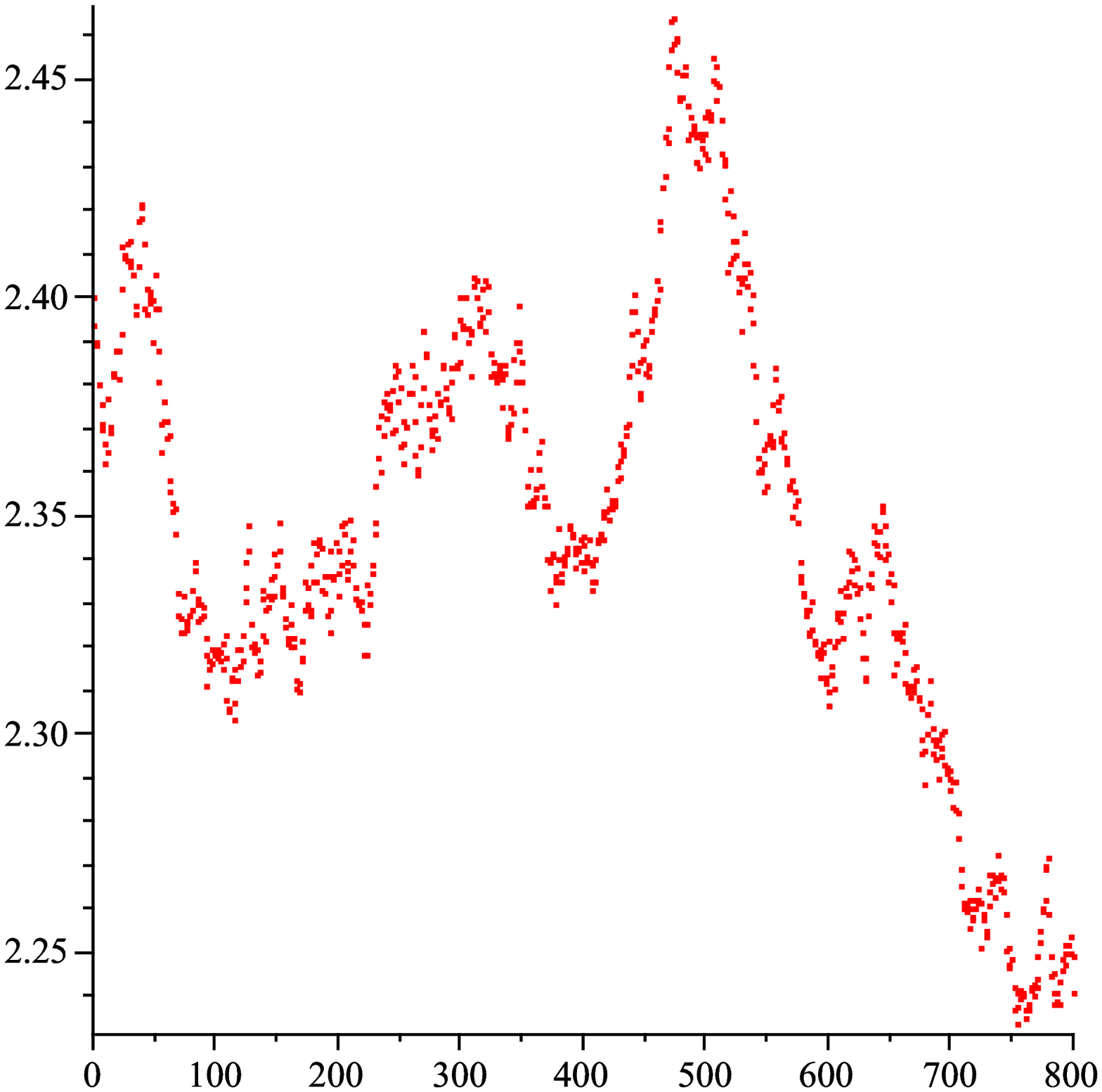} \\
   \footnotesize Figure 17: $(n,y(n))$  in $P_1$  for ODE (\ref{g2})& \footnotesize Figure 18: $(n,y(n,\omega))$ in $P_1$  for SDE (\ref{g3})\\
    &   \\
        &   \\
        \end{tabular}
\end{center}

\begin{center}
\begin{tabular}{cc}
\epsfxsize=6cm \epsfysize=5cm
 \epsffile{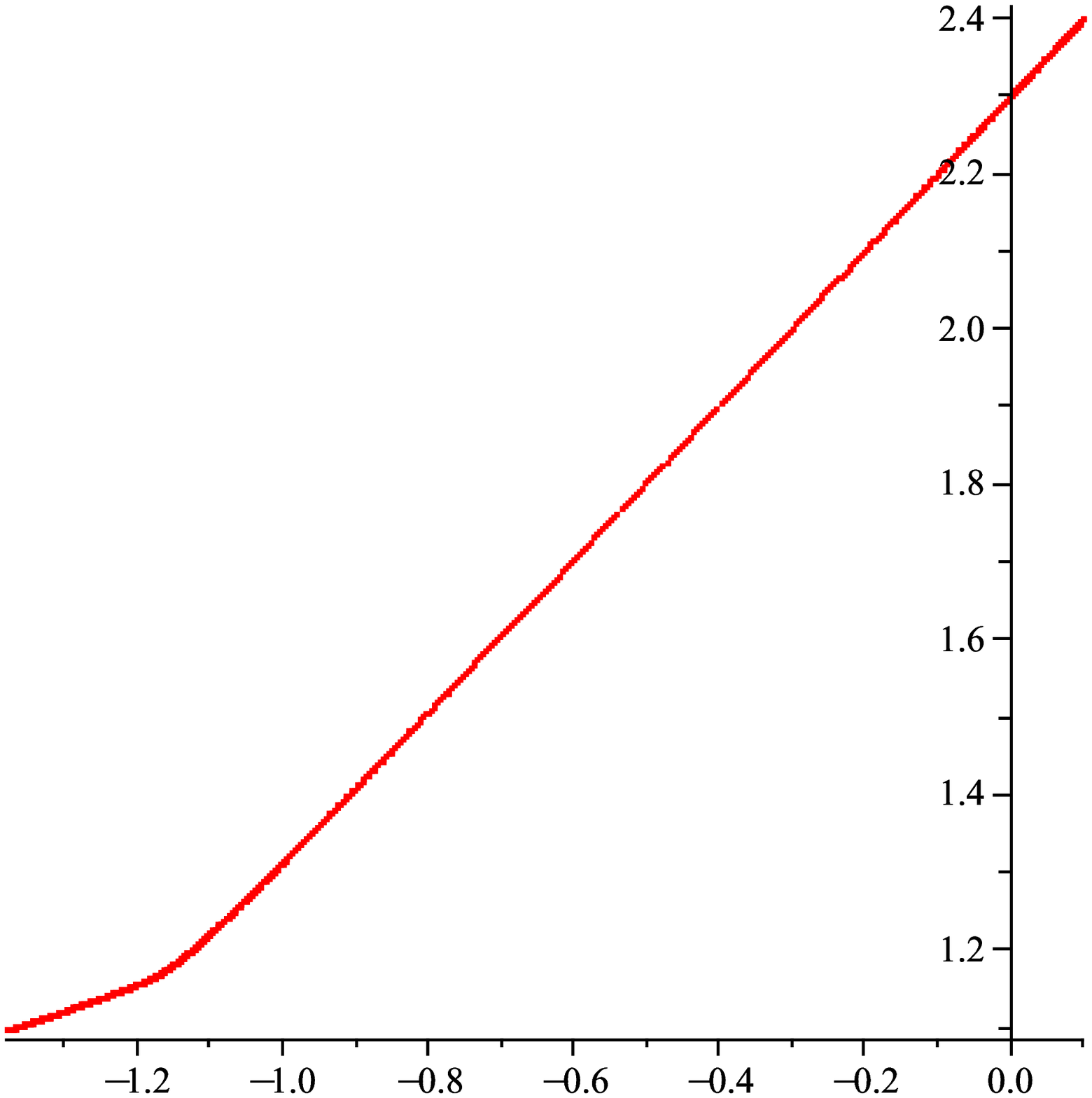}  &
\epsfxsize=6cm \epsfysize=5cm
\epsffile{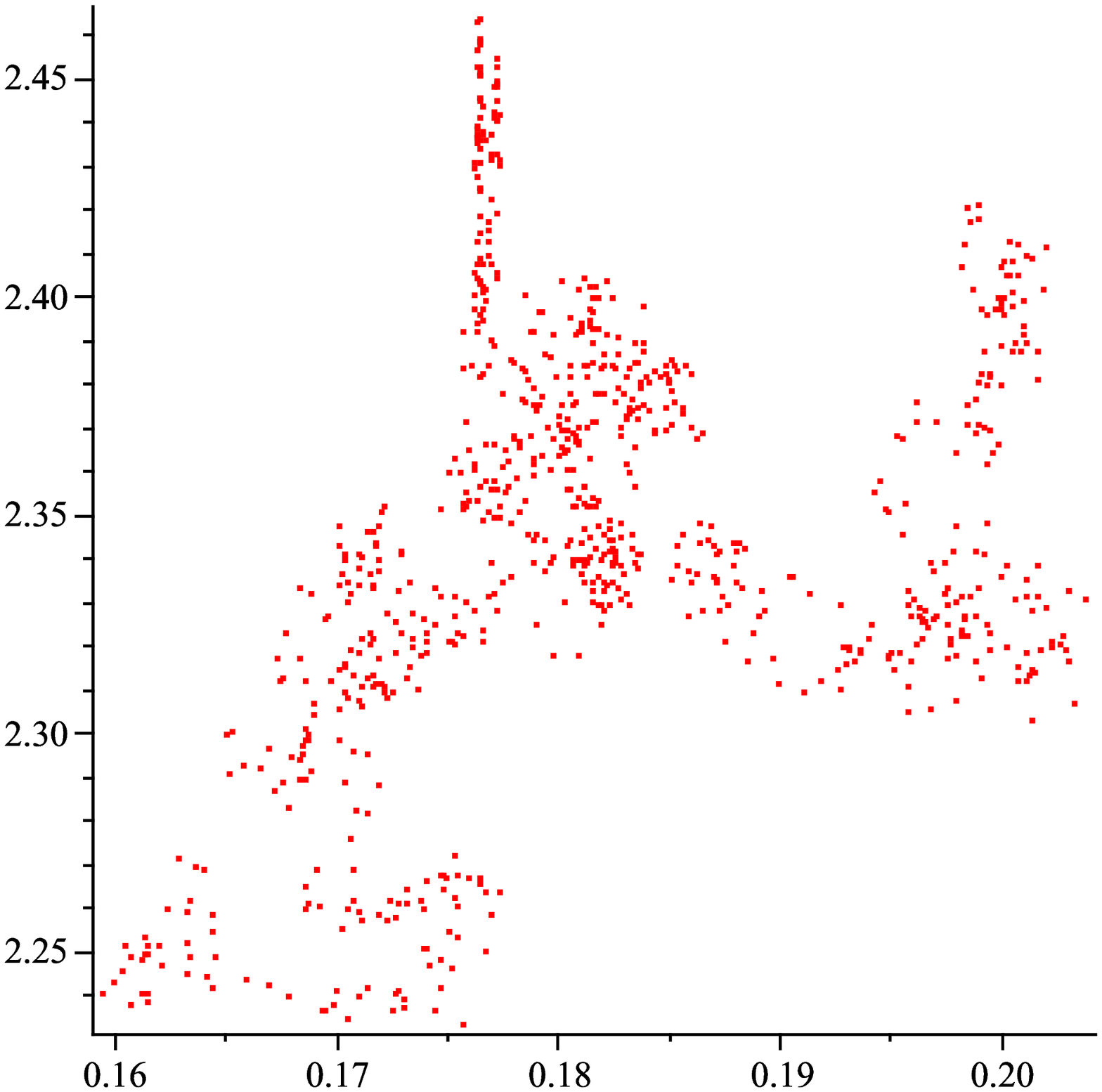} \\
 \footnotesize Figure 19: $(x(n),y(n))$  in $P_1$  for ODE (\ref{g2}) & \footnotesize Figure 20: $(x(n,\omega),y(n,\omega))$ in $P_1$  for SDE (\ref{g3})\\
& \\
& \\
        \end{tabular}
\end{center}

\begin{center}
\begin{tabular}{cc}
 \epsfxsize=6cm \epsfysize=5cm
 \epsffile{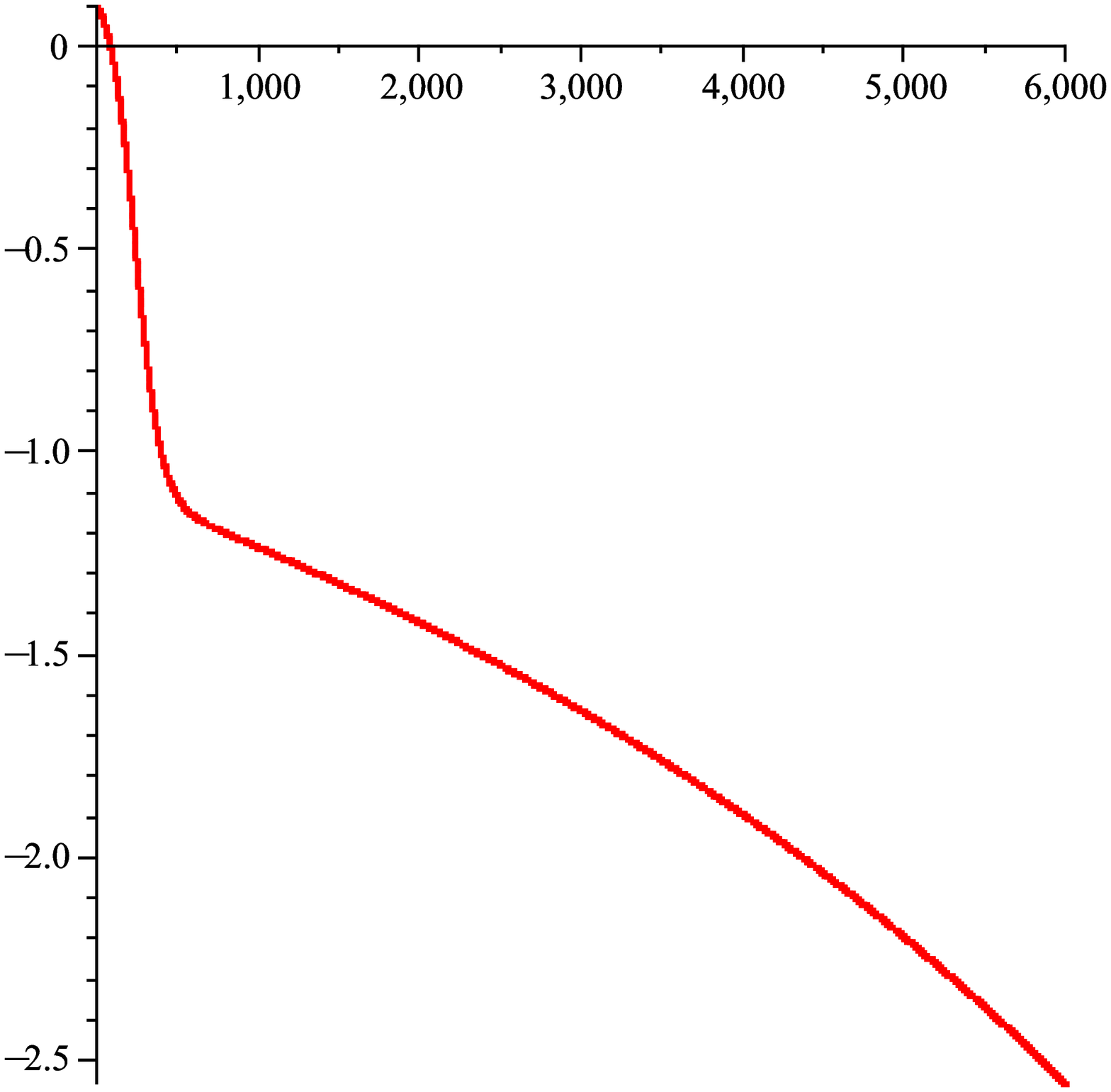}  &
\epsfxsize=6cm \epsfysize=5cm
\epsffile{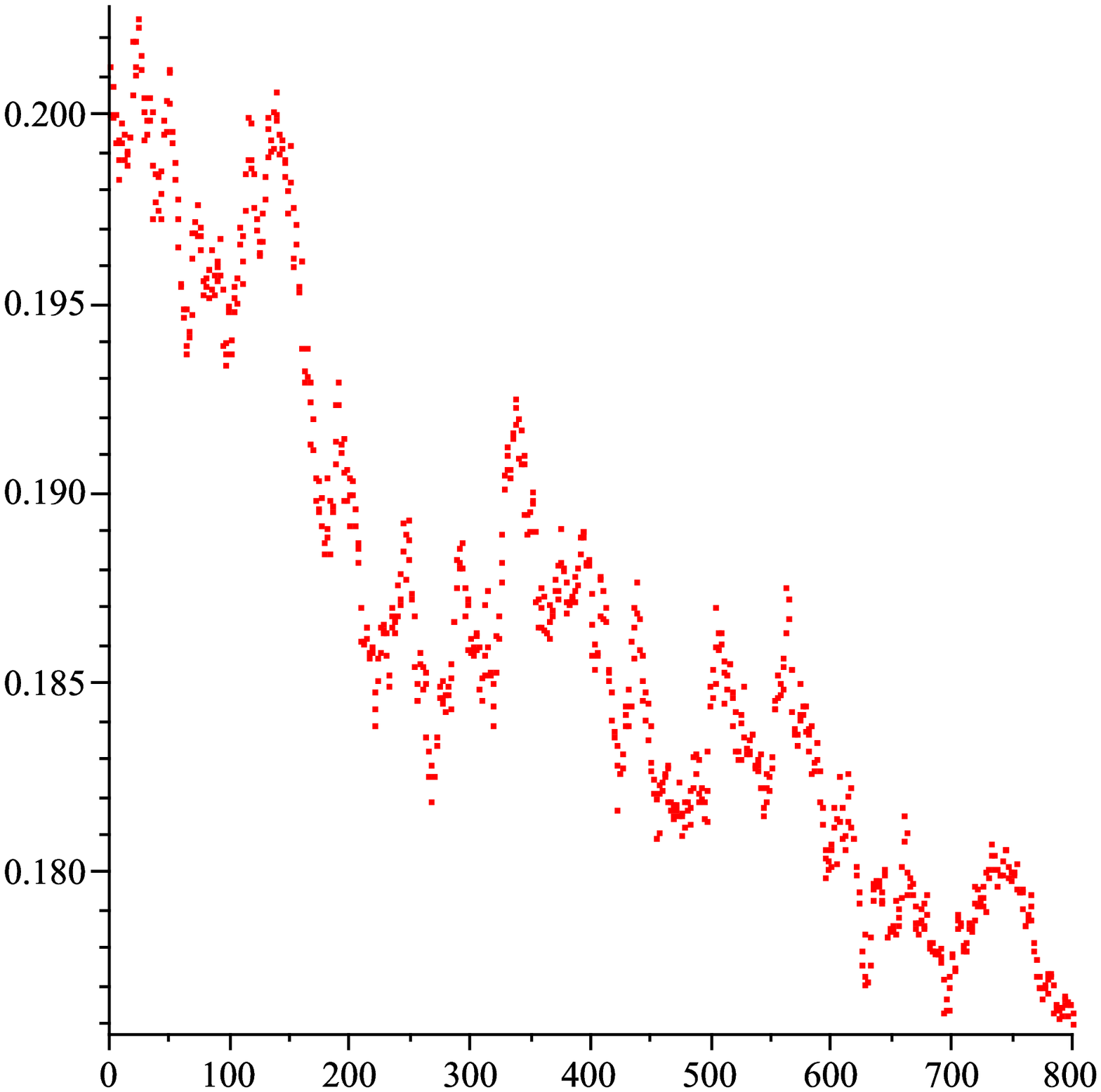} \\
   \footnotesize Figure 21: $(n,x(n))$ in $P_2$  for ODE (\ref{g2})& \footnotesize Figure 22: $(n,x(n,\omega))$ in $P_2$  for SDE (\ref{g3})\\
       &   \\
        &   \\
\end{tabular}
\end{center}

\begin{center}
\begin{tabular}{cc}
\epsfxsize=6cm \epsfysize=5cm
 \epsffile{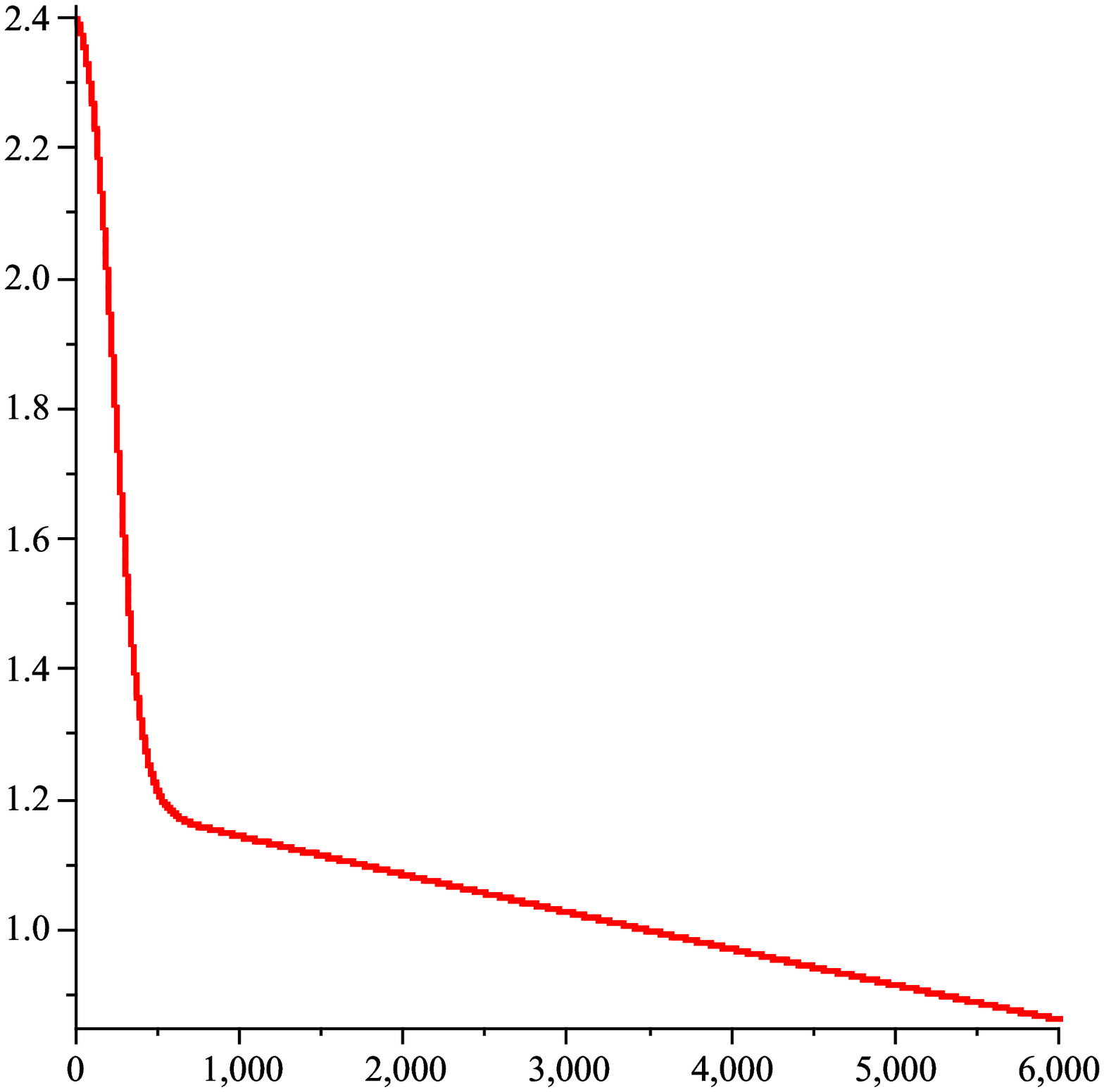}  &
\epsfxsize=6cm \epsfysize=5cm
\epsffile{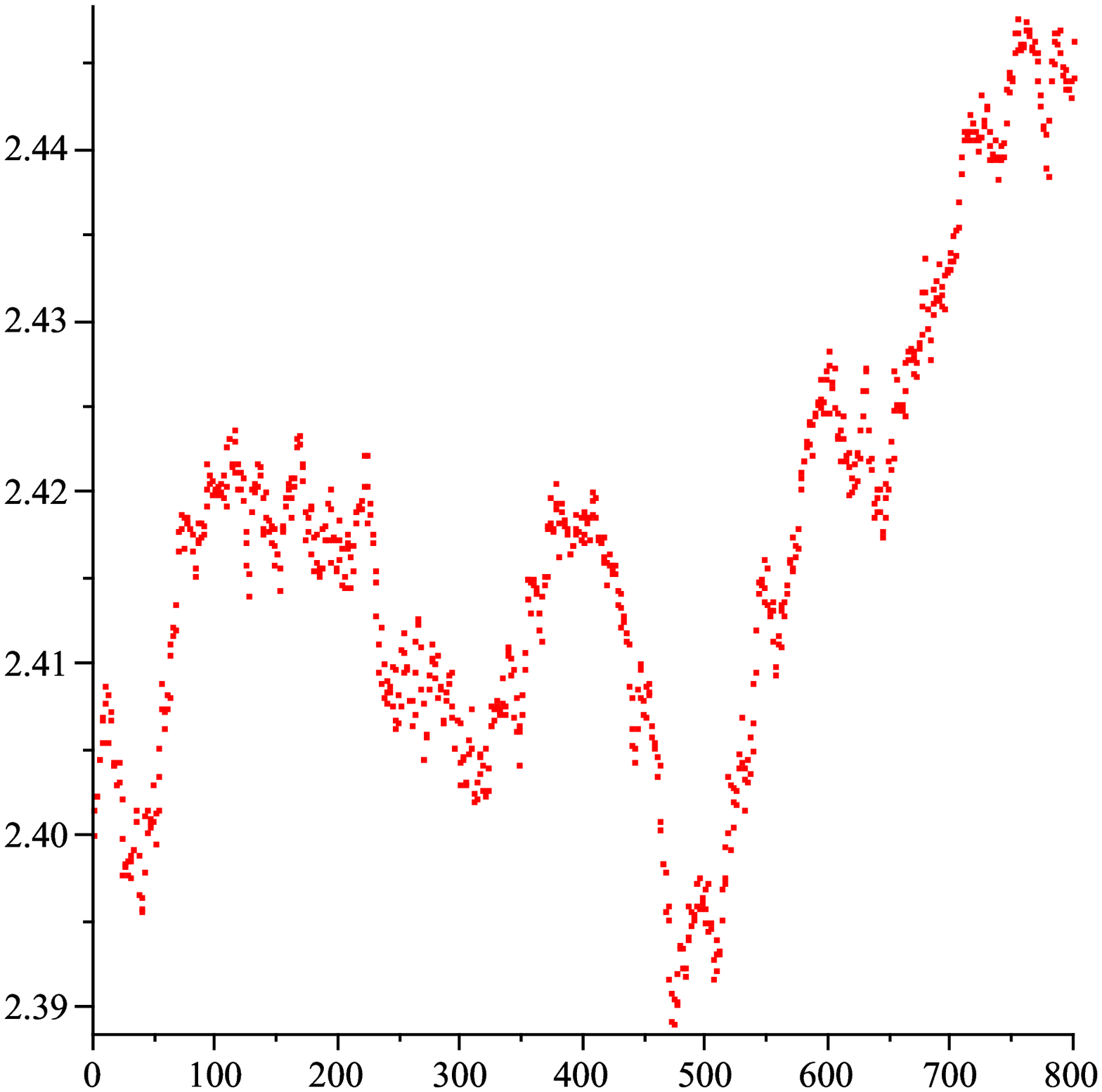} \\
  \footnotesize Figure 23: $(n,y(n))$  in $P_2$  for ODE (\ref{g2})& \footnotesize Figure 24: $(n,y(n,\omega))$ in $P_2$  for SDE (\ref{g3})\\
\end{tabular}
\end{center}

The variation of Lyapunov exponent with the variable parameter
$b_{11}=\alpha$ is given in Figure 23 for $P_1$ and in Figure 24
for $P_2.$

\begin{center}\begin{tabular}{cc}
\epsfxsize=6cm \epsfysize=5cm
 \epsffile{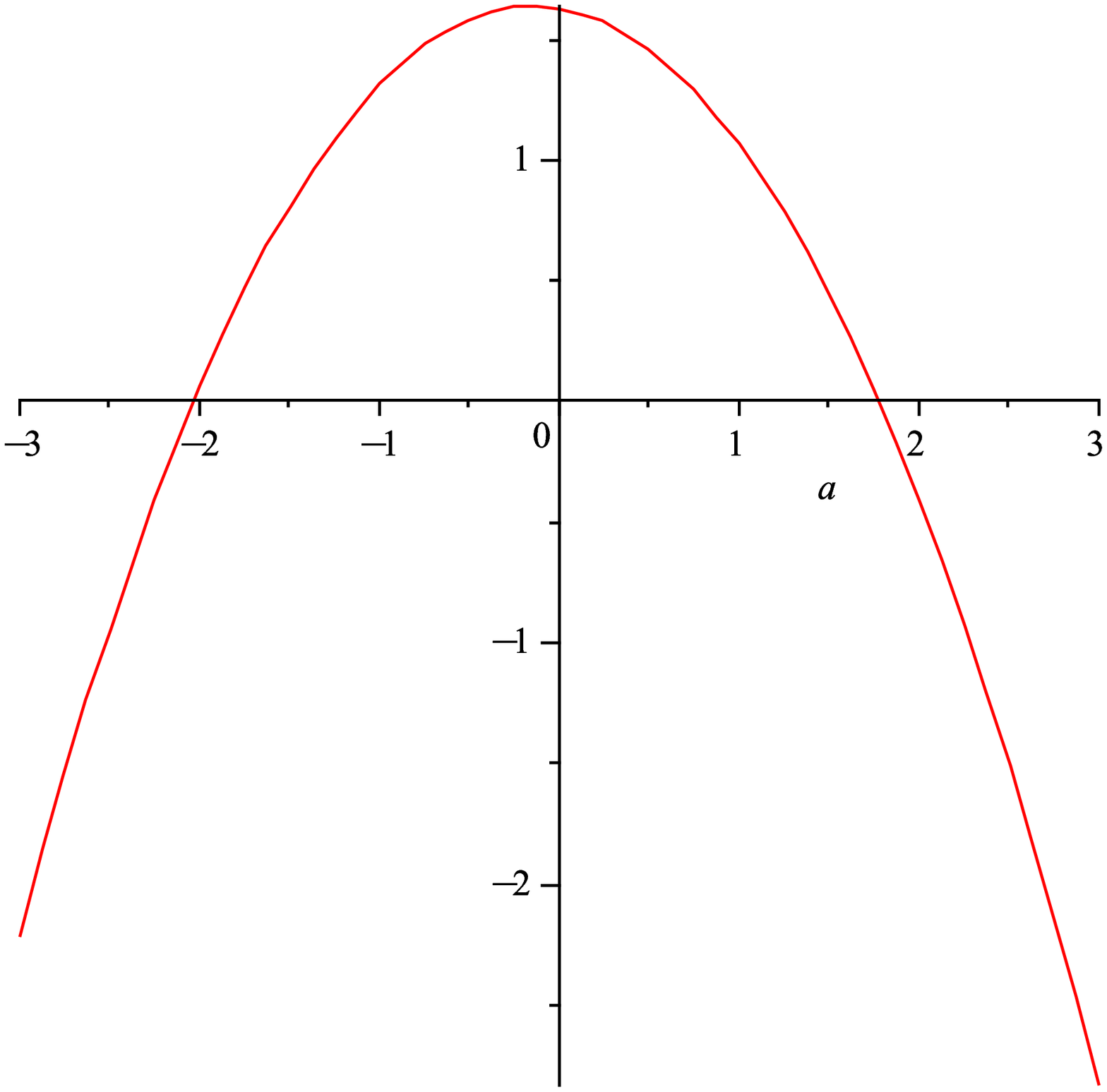}  &
\epsfxsize=6cm \epsfysize=5cm
\epsffile{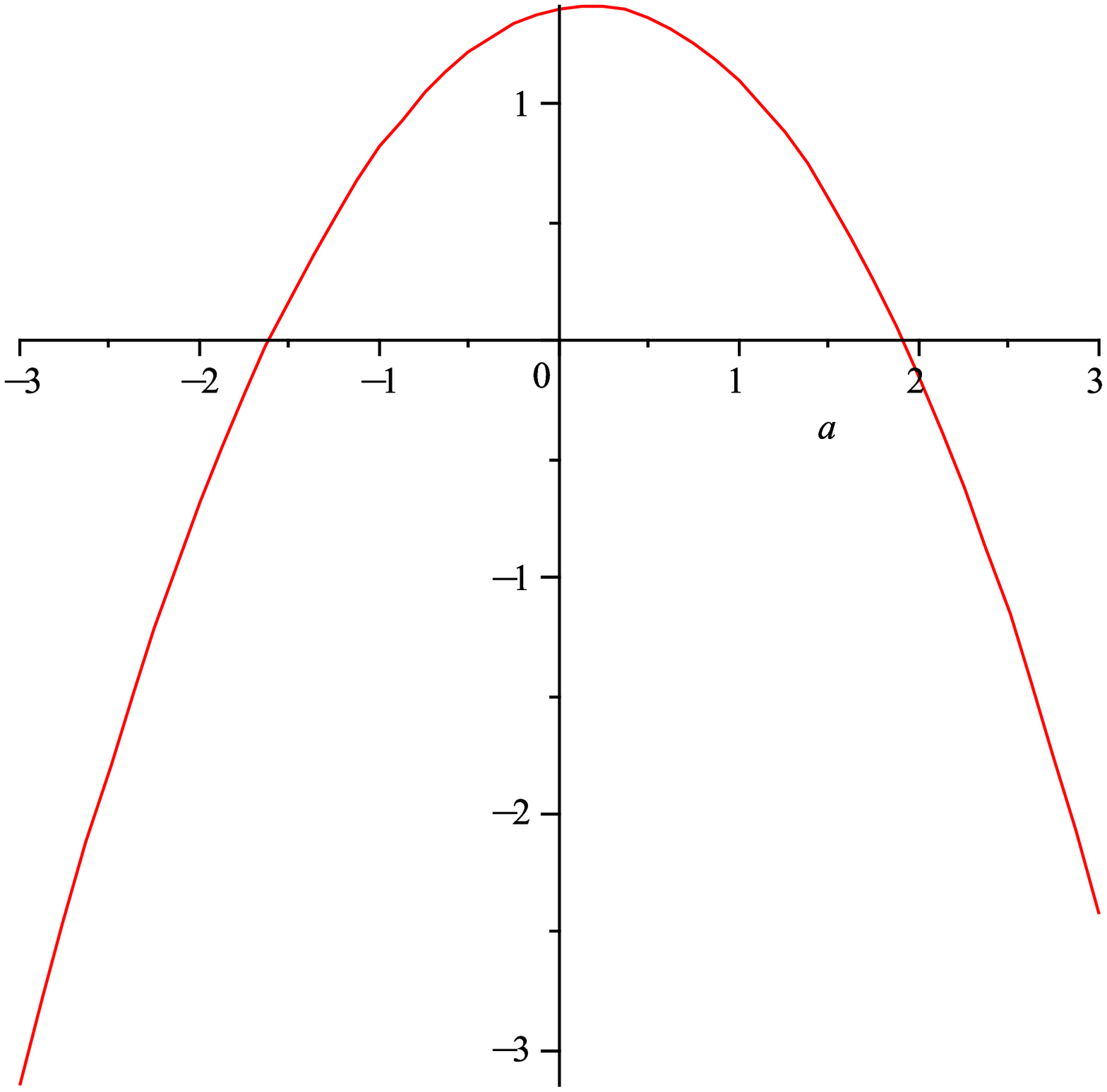} \\
  \footnotesize Figure 25: \small$(\alpha, \lambda(\alpha))$ in $P_1$   for ODE (\ref{g2}) &  \footnotesize Figure 26: $(\alpha, \lambda(\alpha))$ in $P_2$  for SDE (\ref{g3})\\
        &   \\
        &   \\
\end{tabular}
\end{center}

From the figures above, the equilibrium points $P_1$ and $P_2$ are
asymptotically stable for all $\alpha$ such that the Lyapunov
exponents $\lambda(\alpha)<0,$ and unstable otherwise. So, $P_1$
is asymptotically stable for $\alpha\in(-\infty,-2.02)\cup
(1.78,\infty)$ and $P_2$ is asymptotically stable for
$\alpha\in(-\infty,-1.62)\cup(1.88,\infty).$

\subsubsection{Lyapunov function method}

For the system of differential equations that describes Bell
model, the next assertions are true.

\begin{pr}
\begin{description}
    \item[(a)] The matrix of the system of differential equations
    that describes the linearized in $P_2$ is
\end{description}
$$A=\begin{pmatrix}
 0& a_{12} \\
 a_{21} & a_{22}\\
\end{pmatrix},$$ where
$$a_{12}=-\frac{a_2(a_1b_3-a_2b_4)}{a_1b_1-a_2b_2},\quad
a_{21}=\frac{a_1b_1-a_2b_2}{a_2},\quad
a_{22}=-\frac{a_2(b_1b_4+b_2b_3)}{a_1b_1-a_2b_2};$$
    \item[(b)] If $a_1b_1-a_2b_2>0$ and $a_1b_3-a_2b_4>0,$ then
     the equilibrium point $P_2$ is asymptotically stable.
     \end{pr}\hfill $\Box$

The stochastic model is given using a perturbation around the
equilibrium point $P_2(x_2,y_2),$ in the following way
\begin{equation}\label{sistem-stochastic}
\left \{%
\begin{array}{ll}
dx(t)=x(t)(a_1-a_2y(t))dt+\sigma_1(x(t)-x_{P_2})dW^1,\\
dy(t)=[(b_1x(t)-b_3)y(t)-b_2x(t)+b_4]dt+\sigma_2(y(t)-y_{P_2})dW^2,\\
\end{array}
\right.
\end{equation}with $\sigma_1>0, \, \sigma_2>0.$

The linearized of system (\ref{sistem-stochastic}) in $(0,0)$ is
given by
\begin{equation}\label{sistlin}
du(t)=h(u(t))dt+l(u(t))dW(t),
\end{equation} where $u(t)=(u_1(t),u_2(t))^T, \,
W(t)=(W^1(t),W^2(t))^T$ and
\begin{equation}\label{44}
h(u(t))=\begin{pmatrix}
 (a_1-a_2y_2)u_1(t)-x_2a_2u_2(t) \\
 (b_1y_2-b_2)u_1(t)+(b_1x_2-b_3)u_2(t)\\
\end{pmatrix},
\end{equation}
\begin{equation}\label{45}
 l(u(t))=\begin{pmatrix}
 \sigma_1u_1(t)& 0 \\
 0 & \sigma_2u_2(t)\\
\end{pmatrix}.
\end{equation}

We consider the set $D=\{(t\geq 0)\times \mathbb{R}^2\}$ and
$V:D\rightarrow \mathbb{R}$ a function of class $C^1$ with respect
to $t,$ and of class $C^2$ with respect to the other variables. We
study the $p-$exponential stability of the solution $(0,0)$ of the
linearized stochastic system (\ref{sistlin}). Using Theorem
\ref{telin}, from Anexe A1, for the function $V:D \rightarrow
\mathbb{R},$
\begin{equation}\label{v}
V(t,u)=\frac{1}{2}(\omega_1u_1^2+\omega_2u_2^2), \quad
\omega_1,\omega_2\in \mathbb{R}_+,
\end{equation}
we get the following result.

\begin{pr}
If the following relations take place

$$q_1=\omega_1(a_2y_2-a_1-\sigma_1^2)>0, \quad q_2=\omega_2(b_3-b_1y_2-\sigma_2^2)>0,$$
$$b_1y_2-b_2>0, \, \omega_1=\frac{(b_1y_2-b_2)}{a_2y_2}\omega_2,$$
then $$dV(t,u)=-u(t)^TQu(t),$$ with $Q$ given by
$Q=\begin{pmatrix}
 q_1& 0 \\
 0 & q_2\\
\end{pmatrix}.$

The equilibrium point of \emph{(\ref{sistem-stochastic})}
 is asymptotically stable in quadratic square ($p=2$).
\end{pr}
\textbf{Proof:} From (\ref{44}), (\ref{45}) and (\ref{v}), we
get
\begin{eqnarray*}
dV(t,u)&=&\begin{pmatrix}
 (a_1-a_2y_2)u_1(t)-x_2a_2u_2 \\
 (b_1y_2-b_2)u_1+(b_1x_2-b_3)u_2\\
\end{pmatrix}^T
\begin{pmatrix}
 \omega_1u_1 \\
 \omega_2u_2\\
\end{pmatrix}+\begin{pmatrix}
 \sigma_1^2\omega_1u_1^2& 0 \\
 0 & \sigma_2^2\omega_2u_2^2\\
\end{pmatrix}\\
&=&-q_1u_1^2-q_2u_2^2+(\omega_2(b_1y_2-b_2)-\omega_1x_2a_2)u_1u_2.
\end{eqnarray*}
If the relations from the proposition take place, then we get
$$dV(t,u)=-u(t)^TQu(t).$$

The matrix $Q$ is symmetric and positive defined  and has positive
eigenvalues $r_1=q_1$ and $r_2=q_2.$  Let $q_{m}$ be
$q_m=min\{q_1,q_2\}.$ Results that
$$LV(t,u)\leq -q_m\|u(t)\|^2.\,$$ and so the equilibrium point is asymptotically
stable in square mean.

\hfill $\Box$

Let us choose the same parameters values for $a_1, \, a_2,\, b_1, \, b_2, \, b_3$ as on the simulation of Lyapunov exponents. We use Maple 12 software for the implementation of the second order Euler method. We observe from the following graphics that the solution
trajectories represent the stable characteristic, which validate our theoretical discussion for the  system of
differential equation (\ref{sistem-stochastic}), for the equilibrium point $P_2.$

\begin{center}\begin{tabular}{cc}
\epsfxsize=6cm \epsfysize=5cm
 \epsffile{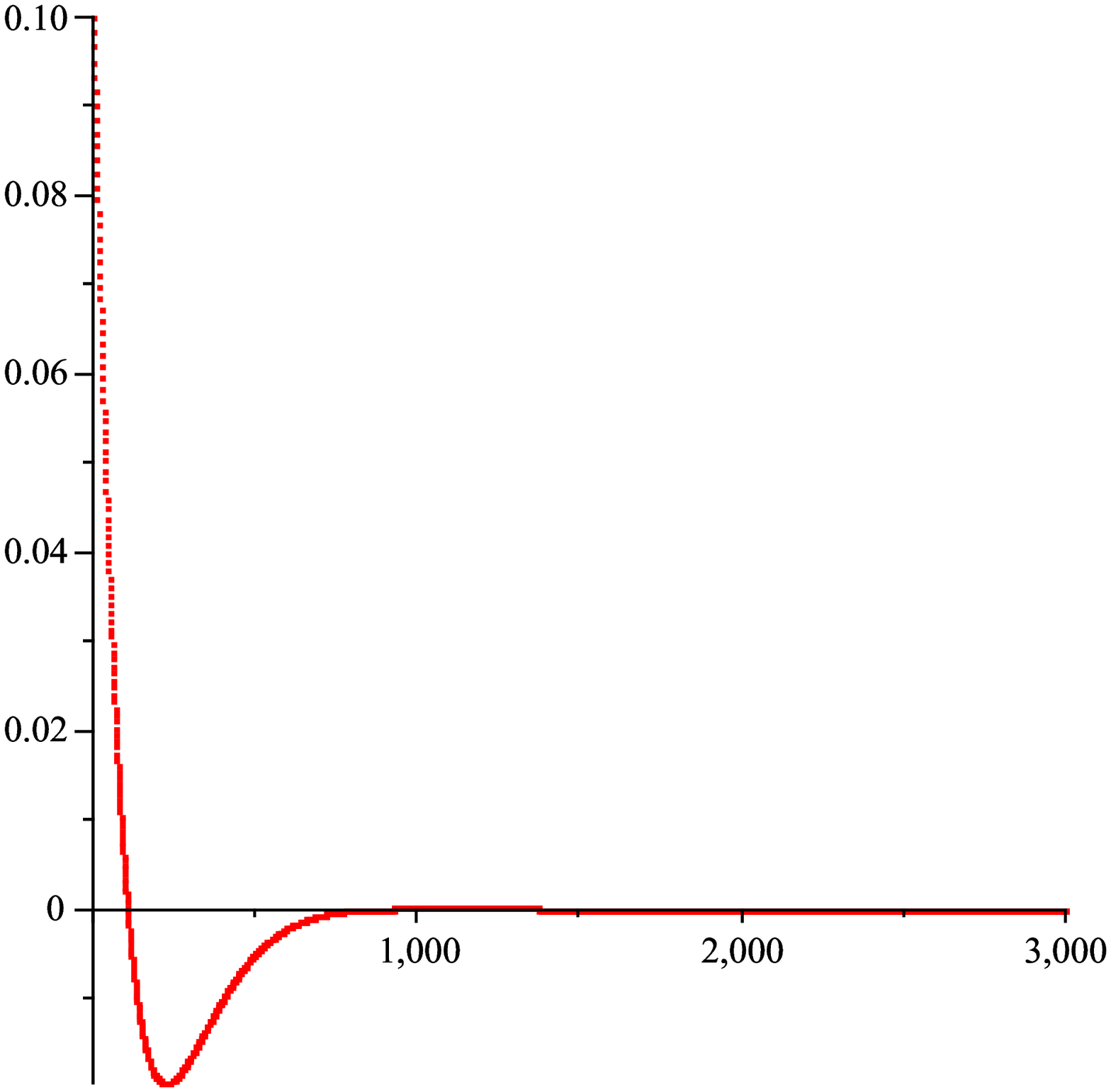}  &
\epsfxsize=6cm \epsfysize=5cm
\epsffile{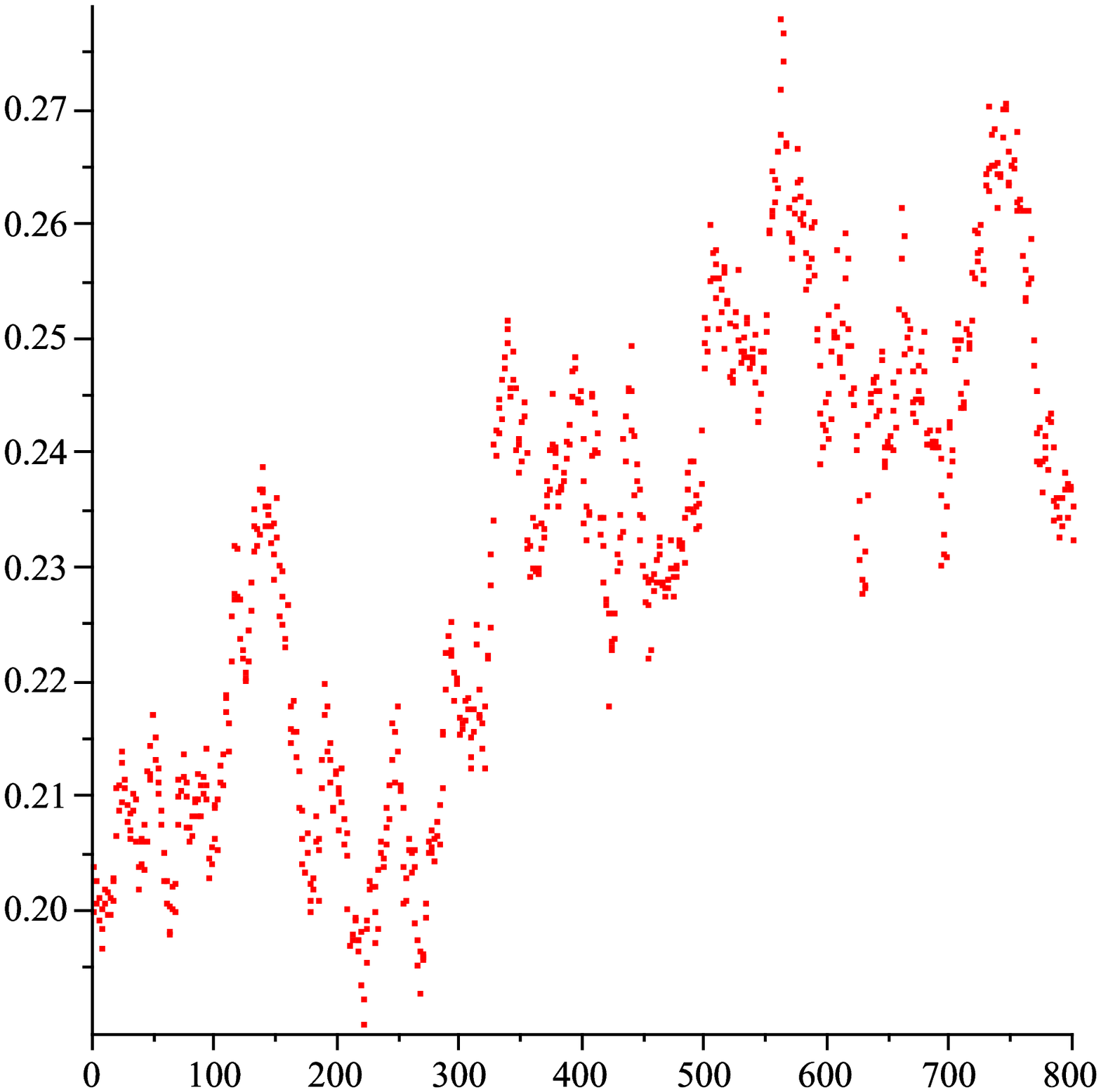} \\
 \footnotesize Figure 27: $(n,x(n))$  in $P_1$ for ODE (\ref{g2}) &\footnotesize Figure 28: $(n,x(n,\omega))$ in $P_1$  for SDE (\ref{sistlin})\\
 &   \\
        &   \\
        \end{tabular}
\end{center}

\begin{center}\begin{tabular}{cc}
\epsfxsize=6cm \epsfysize=5cm
 \epsffile{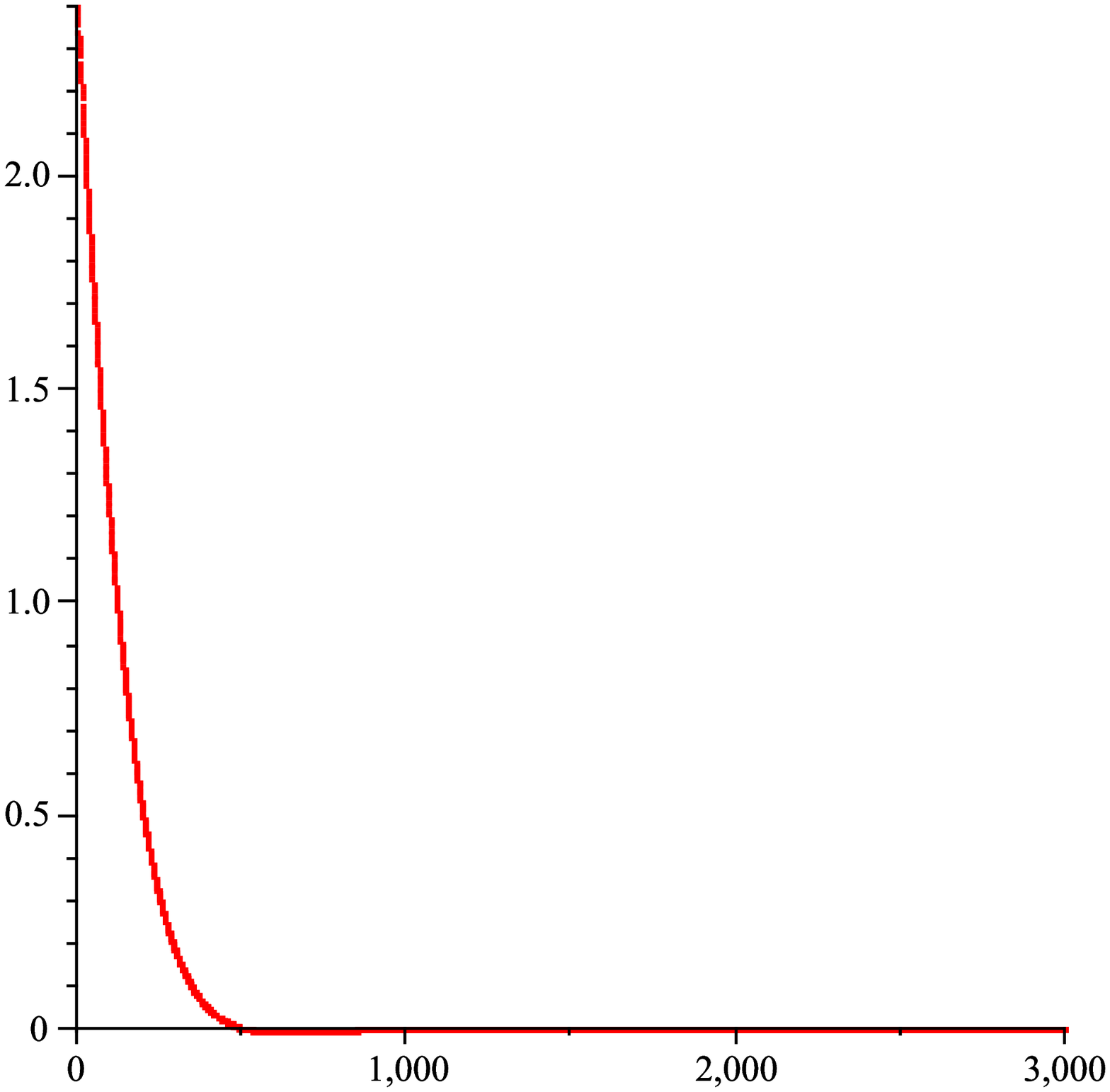}  &
\epsfxsize=6cm \epsfysize=5cm
\epsffile{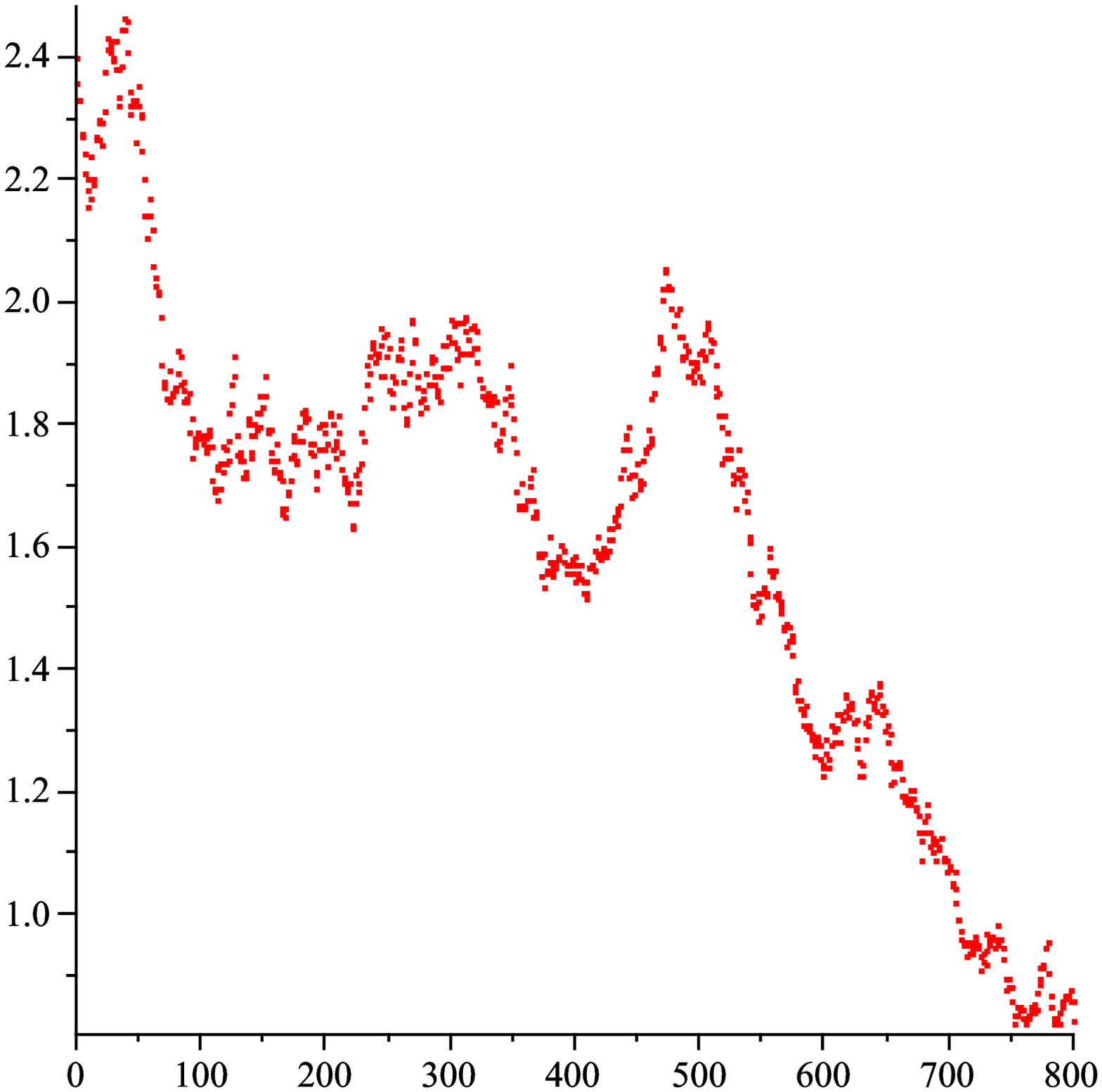} \\
 \footnotesize Figure 29: $(n,x(n))$  in $P_1$ for ODE (\ref{g2}) &\footnotesize Figure 30: $(n,x(n,\omega))$ in $P_1$  for SDE (\ref{sistlin})\\
 &   \\
        &   \\
        \end{tabular}
\end{center}

\section{Conclusions}

In this paper we focused on two important tumor-immune systems,
presented from stochastic point of view: a Kuznetsov-Taylor model
and Bell model, that belongs to a general family of tumor-immune
stochastic systems. We have determined the equilibrium points and
we have calculated the Lyapunov exponents. A computable algorithm
is presented in Annexe A1. These exponents help us to decide
whether the stochastic model is stable or not. We have also proved stochastic 
stability by constructing a proper Lyapunov function, under a well chosen conditions.
All our results were also proved using graphical implementation.
For numerical
simulations we have used the second order Euler scheme presented in detail in
Annexe A2 and the implementation of this algorithm was done in
Maple 12.

In a similar way can be studied other models that derive from model given by (\ref{g2}). The model given by the SDE (\ref{g3}) allows the control of the system given by the ODE (\ref{3}), with a stochastic process.

\section*{Annexe}

\subsection*{A1  Lyapunov exponents and stability in stochastic 2-
 dimensional structures.}

\subsubsection*{Lyapunov exponent method}

\par The behavior of a deterministic dynamical system which is
disturbed by noise may be modeled by a stochastic differential
equation (SDE). In many practical situations, perturbations are
generated by wind, rough surfaces or turbulent layers are
expressed in terms of white noise, modeled by Brownian motion. The
stochastic stability has been introduced in \cite{Jed} and is
characterized by the negativeness of Lyapunov exponents.  But it
is not possible to determine this exponents explicitly. Many
numerical approaches have been proposed, which generally used
simulations of stochastic trajectories.

Let $(\Omega, \mathcal{F},\mathcal{P})$ a probability space. It is
assumed that the $\sigma-$algebra $(\mathcal{F}_{t})_{t\geq 0}$ such
that
$$\mathcal{F}_s\subset\mathcal{F}_t\subset\mathcal{F}, \, \forall \, s\leq t, \,
s,t\in I,$$ where $I=[0,T], \, T\in (0,\infty).$

Let $\{x(t,\omega)=(x_1(t),x_2(t))\}_{t\geq 0}$ be a stochastic
process, solution of the system of It\^o differential equations,
formally written as
\begin{equation}\label{aa1}
dx_i(t,\omega)=f_i(x(t,\omega))dt+g_i(x(t,\omega))dW(t,\omega), \,
i=1,2,
\end{equation}
with initial condition $x(0)=x_0$ is interpreted in the sense that
\begin{equation}\label{aa2}
x_i(t,\omega)=x_{i0}(t,\omega)+\int_0^tf_i(x(s,\omega))ds+\int_0^t
g_i(x(s,\omega))dW(s,\omega), \, i=1,2,
\end{equation}
for almost all $\omega\in\Omega$ and for each $t>0,$ where
$f_i(x)$ is a drift function, $g_i(x)$ is a diffusion function,
$\int_0^t f_i(x(s))ds, \, i=1,2$ is a Riemann integral and
$\int_0^t g_i(x(s))dW(s), \, i=1,2$ is an It\^o integral. It is
assumed that $f_i$ and $g_i, \, i=1,2$ satisfy thce conditions of
existence of solutions for this SDE with initial conditions
$x(0)=a_0\in \mathbb{R}^n.$

Let $x_0=(x_{01},x_{02})\in  \mathbb{R}^2$ be a solution of the
system
\begin{equation}\label{aa3}
f_i(x_0)=0, \, i=1,2.
\end{equation}
The functions $g_i$ are chosen such that
$$g_i(x_0)=0, \, i=1,2.$$
In the following, we will consider
\begin{equation}\label{aa4}
g_i(x)=\sum_{j=1}^2 b_{ij}(x_j-x_{j0}), \, i=1,2,
\end{equation}where $b_{ij}\in \mathbb{R}, \, i,j=1,2.$

The linearized of the system (\ref{aa2}) in $x_0$ is given by
\begin{equation}\label{aa5}
du(t)=Au(t)dt+ Bu(t)dW(t),
\end{equation}where
\begin{equation}\label{aa6}
u(t)=\begin{bmatrix}
 u_1(t) \\
 u_2(t) \\
\end{bmatrix}, \quad A=\begin{bmatrix}
 a_{11}  & a_{12} \\
 a_{21}  & a_{22}  \\
\end{bmatrix}, \quad B=\begin{bmatrix}
 b_{11}  & b_{12} \\
 b_{21}  & b_{22}  \\
\end{bmatrix}
\end{equation}
\begin{equation}\label{aa7}
a_{ij}=\frac{\partial f_i}{\partial x_j}\Big|_{x=x_0},\quad
b_{ij}=\frac{\partial g_i}{\partial x_j}\Big|_{x=x_0}, \quad i,j=1,2.
\end{equation}

The Oseledec multiplicative ergodic theorem \cite{Ose} asserts the
existence of two non-random Lyapunov exponents $\lambda_2\leq
\lambda_1=\lambda.$  The top Lyapunov exponent is given by
\begin{equation}\label{aa8}
\lambda=\mathop{\lim}
\limits_{t\to\infty}\sup \frac{1}{t}\log\sqrt{u_1(t)^2+u_2(t)^2}.
\end{equation}

By applying the change of coordinates
$$u_1(t)=r(t)\cos \theta(t), \, u_2(t)=r(t)\sin \theta(t),$$
 for (\ref{aa5}) and by using the It\^o formula for
 \begin{eqnarray*}
h_1(u_1,u_2)&=&\frac{1}{2}\log(u_1^2+u_2^2)=\log (r),\\
h_2(u_1,u_2)&=&\arctan\Big(\frac{u_2}{u_1}\Big)=\theta,
\end{eqnarray*}
 result
the stochastic equations written in the integral form.
\begin{equation}\label{aap1}
\log\Big(\frac{r(t)}{r(0)}\Big)=\int_0^t
q_1(\theta(s))+\frac{1}{2}(q_4(\theta(s))^2-q_2(\theta(s))^2)ds+
\int_0^tq_2(\theta(s))dW(s),
\end{equation}
\begin{equation}\label{aap2}
\theta(t)=\theta(0)+\int_0^t
(q_3(\theta(s))-q_2(\theta(s))q_4(\theta(s)))ds+
\int_0^tq_4(\theta(s))dW(s),
\end{equation}where
\begin{equation}\label{aap3}
\begin{array}{ll}
q_1(\theta)=a_{11}\cos^2\theta+(a_{12}+a_{21})\cos\theta\sin\theta+a_{22}\sin^2\theta,\\
q_2(\theta)=b_{11}\cos^2\theta+(b_{12}+b_{21})\cos\theta\sin\theta+b_{22}\sin^2\theta,\\
q_3(\theta)=a_{21}\cos^2\theta+(a_{22}-a_{11})\cos\theta\sin\theta-a_{12}\sin^2\theta,\\
q_4(\theta)=b_{21}\cos^2\theta+(b_{22}-b_{11})\cos\theta\sin\theta-b_{12}\sin^2\theta.
\end{array}
\end{equation}
As the expectation of the It\^o stochastic integral is null,
$$E\Big(\int_0^t q_2(\theta(s))dW(s)\Big)=0,$$ the Lyapunov exponent is
given by
$$\lambda=\mathop{\lim}\limits_{t\to\infty}\frac{1}{t}\log\Big(\frac{r(t)}{r(0)}\Big)=
\mathop{\lim}\limits_{t\to\infty}\frac{1}{t}E\Big(\int_0^t
(q_1(\theta(s))+\frac{1}{2}(q_4(\theta(s))^2-q_2(\theta(s))))ds\Big).$$
Applying the Oseledec theorem, if $r(t)$ is ergodic, results that
\begin{equation}\label{aap4}
\lambda=\int_0^t(q_1(\theta)+\frac{1}{2}(q_4(\theta)^2-q_2(\theta)))p(\theta)d\theta,
\end{equation}where $p(\theta)$ is the probability density of
the process $\theta.$

The probability density is the solution $p(t,\theta)$  of Fokker-Planck equation associated to equation (\ref{aap2})
\begin{equation}\label{aap5}
\frac{\partial p(t,\theta)}{\partial t}+\frac{\partial }{\partial
\theta}(q_3(\theta)-q_2(\theta)q_4(\theta)p(t,\theta))-\frac{1}{2}\frac{\partial^2
}{\partial \theta^2}(q_4(\theta)^2p(t,\theta))=0.
\end{equation}

If $p(t,\theta)=p(\theta),$ then the stationary solution of  (\ref{aap5}) satisfies the first order differential equation
\begin{equation}\label{aap6}
(-q_3(\theta)+q_1(\theta)q_4(\theta)+q_2(\theta)q_5(\theta))p(\theta)+\frac{1}{2}
q_4(\theta)^2\dot{p} (\theta)=p_0,
\end{equation}where $\dot{p} (\theta)=\frac{dp}{d\theta}$ and
\begin{equation}\label{aap7}
q_5(\theta)=-(b_{12}+b_{21})\sin (2\theta)-(b_{22}-b_{11})\cos
(2\theta).
\end{equation}

\begin{pr}
If $q_4(\theta)\neq 0,$ the solution of equation
\emph{(\ref{aap6})} is given by
\begin{equation}\label{aap8}
p(\theta)=\frac{K}{D(\theta)q_4(\theta)^2}(1+\eta \int_0^\theta
D(u)du),
\end{equation}
where $K$ is determined by the normality condition
\begin{equation}\label{aap9}
\int_0^{2\pi} p(\theta)d\theta=1,
\end{equation}and
\begin{equation}\label{aap10}
\eta=\frac{D(2\pi)-1}{\int_0^{2\pi} D(u)du}.
\end{equation}The function $D$ is given by
\begin{equation}\label{aap11}
D(\theta)=exp(-2
\int_0^\theta\frac{q_3(u)-q_2(u)q_4(u)-q_4(u)q_5(u)}{q_4(u)^2}du).
\end{equation}\hfill $\Box$
\end{pr}

A numerical solution of the phase distribution could be performed
by a simple backward difference scheme. The function $p(\theta)$ can be determined numerically by using the following algorithm.

Let us consider $N\in \mathbb{R}_+,$  $h=\frac{\pi}{N}$ and
\begin{equation}\label{aap12}
\begin{array}{ll}
q_1(i)=a_{11}\cos^2(ih)+(a_{12}+a_{21})\cos(ih)\sin(ih)+a_{22}\sin^2(ih),\\
q_2(i)=b_{11}\cos^2(ih)+(b_{12}+b_{21})\cos(ih)\sin(ih)+b_{22}\sin^2(ih),\\
q_3(i)=a_{21}\cos^2(ih)+(a_{22}-a_{11})\cos(ih)\sin(ih)-a_{22}\sin^2(ih),\\
q_4(i)=b_{21}\cos^2(ih)+(b_{22}-b_{11})\cos(ih)\sin(ih)-b_{12}\sin^2(ih),\\
q_5(i)=-(b_{12}+b_{21})\sin (2ih)-(b_{22}-b_{11})\cos (2ih), \quad i=0,1,...,N.
\end{array}
\end{equation}
The sequence $(p(i))_{i=0,...,N}$ is given by
$$p(i)=(p(0)+\frac{q_4(i)^2p(i-1)}{2h})F(i),$$ where
$$F(i)=\frac{2h}{2h(-q_3(i)+q_2(i)q_4(i)+q_4(i)q_5(i))+q_4(i)^2}.$$
The Lyapunov exponent is $\lambda=\lambda(N),$ where
$$\lambda(N)=\sum_{i=1}^N(q_1(i)+\frac{1}{2}(q_4(i)^2-q_2(i)^2))p(i)h.$$

From (\ref{aap3}) and (\ref{aap8}) we get the following proposition.

\begin{pr}
If the coefficients of the the matrix $B$ are given by
$$b_{11}=\alpha, \, b_{12}=-\beta, \, b_{21}=\beta, \,
b_{22}=\alpha,$$ then
the Lyapunov exponent is given by
$$\lambda=\frac{1}{2}(a_{11}+a_{22}+\beta^2-\alpha^2)+\frac{1}{2}(a_{11}-a_{22})D_2+
\frac{1}{2}(a_{21}+a_{12})E_2,$$ where
$$
D_2=\int\limits_{0}^{2\pi}\cos(2\theta)p(\theta)d\theta,\quad
E_2=\int\limits_{0}^{2\pi}\sin(2\theta)p(\theta)d\theta.
$$

$$p(\theta)=Kg(\theta), \quad K=\frac{1}{\int_0^{2\pi} g(\theta)d\theta}, \,$$ $$g(\theta)=\frac{1}{\beta^2}exp(\frac{1}{\beta^2}((a_{21}-a_{12}-\alpha\beta)\theta+
\frac{1}{2}(a_{11}-a_{22})\cos
(2\theta)+\frac{1}{2}(a_{21}-a_{12})\sin (2\theta))).$$
 \hfill $\Box$
\end{pr}

\subsubsection*{Lyapunov function method}

Let us consider the stochastic system of differential equations given by

\begin{equation}\label{a121}
dx_i(t)=f_i(x(t))dt+g_i(x(t))dW_i(t), \quad i=1,2,
\end{equation}
where $W_1, \, W_2$ are Wiener processes. Let $D=(0,\infty)\times
\mathbb{R}^2,$ and $V:D\rightarrow \mathbb{R}$ a continuous
function with respect to the first component and of the class
$C^2$ with respect to the other components. Let consider the
differential operator given by
\begin{equation}\label{a122}
LV(t,x)=\frac{\partial V(t,x)}{\partial t}+\sum_{i=1}^2 f_i(x)\frac{\partial V(t,x)}{\partial x_i}+\frac{1}{2}\sum_{i=1}^2 \sum_{j=1}^2 g_i(x)g_j(x)\frac{\partial^2 V(t,x)}{\partial x_i \partial x_j}.
\end{equation}

We suppose that $x_0=0$ is an equilibrium point for (\ref{a121}), that means
\begin{equation}\label{a123}
f_i(0)=g_i(0)=0, \quad i=1,2.
\end{equation}

The theorem that gives the conditions for stability of the trivial
solution $x_0=0$ in terms of Lyapunov function is given in
\cite{Schu}.

\begin{te}
If there is a function $V:U\rightarrow \mathbb{R}$ and two continuous functions  $u,v:\mathbb{R}_+ \rightarrow \mathbb{R}_+$ and $k>0,$ such that for each $\|x\|<k,$ we have
\begin{equation}
u(\|x\|)<V(x,t)<v(\|x\|),
\end{equation} then
\begin{description}
    \item[(i)] If $LV(t,x)\leq 0, \, x\in (0,k),$ then the solution $x_0=0$ of \emph{(\ref{a121})} is stable in probability,
    \item[(ii)] If there is a continuous function $c:\mathbb{R}_+ \rightarrow \mathbb{R}_+$ such that $$LV(t,x)\leq -c(\|x\|),$$ then the solution $x_0=0$ of \emph{(\ref{a121})} is asymptotically  stable. \hfill $\Box$
\end{description}
\end{te}

In general, the functions $f_i, \, g_{i}, \, i=1,2,$ are nonlinear
and the above theorem is hard to use. That is why we use the
linearization method for the system (\ref{a121}) around the
equilibrium point.

The linearized system of stochastic differential equation of (\ref{a121}) is given by

\begin{equation}\label{a124}
\left \{%
\begin{array}{ll}
du_1(t)=(a_{11}u_1(t)+a_{12}u_2(t))dt+(b_{11}u_1(t)+b_{12}u_2(t))dW_1, \\
du_2(t)=(a_{21}u_1(t)+a_{22}u_2(t))dt+(b_{21}u_1(t)+b_{22}u_2(t))dW_2. \\
\end{array}
\right.
\end{equation}

We consider $D=\{(t\geq0) \times \mathbb{R}^2\}$ and $V:D\rightarrow \mathbb{R}$ a continuous function with respect to $t$ and of the class $C^2$ with respect to the other components. The theorem that gives the condition that the trivial solution of (\ref{a124}) is exponential $p-$stable is given in \cite{p-sta}.

\begin{te}\label{telin}
If the function $V:D\rightarrow \mathbb{R}$ satisfies the inequalities
$$k_1\|u\|^p\leq V(t,u)\leq k_2 \|u\|^p,$$
$$LV(t,u)\leq -k_3\|u\|, \quad k_i>0, \, p>0,$$
then the trivial solution of \emph{(\ref{a124})} is exponentially $p-$stable for $t\geq 0.$ \hfill $\Box$
\end{te}

In concrete problems, the next theorem is used.

\begin{te}\label{telin-1}
If the function $V:D\rightarrow \mathbb{R}$ satisfies
\begin{description}
    \item[(i)] $LV(u)\leq 0,$ then the trivial solution is stable in probability;
    \item[(ii)] $LV(u)\leq -c(\|u\|),$ where $c:\mathbb{R}_+\rightarrow \mathbb{R}_+$ is a continuous function, then the trivial solution is asymptotically stable;
        \item[(iii)] $LV(u)\leq -q^TQ q,$ where $Q$ is a symmetric matrix, positive defined, then the trivial solution is stable in mean square value.
\end{description}\hfill $\Box$
\end{te}

For (\ref{a124}), the expression of the differential operator  $LV$ is given by

\begin{eqnarray*}
LV(t,u)&=&(a_{11}u_1+a_{12}u_2)\frac{\partial V(t,u)}{\partial u_1}+
(a_{21}u_1+a_{22}u_2)\frac{\partial V(t,u)}{\partial u_2}\\
&+&\frac{1}{2}\Big[(b_{11}u_1+b_{12}u_2)^2\frac{\partial^2 V(t,u)}{\partial u_1^2}+(b_{21}u_1+b_{22}u_2)^2\frac{\partial^2 V(t,u)}{\partial u_2^2}\Big].
\end{eqnarray*}

\subsection*{A2 The Euler scheme.}

In general 2-dimensional case, the Euler scheme has the form:

\begin{equation}\label{aaa1}
x_i(n+1)=x_i(n)+f_i(x(n))h+g_i(x(n))G_i(n), \, i=1,2,
\end{equation}with Wiener process increment
$$G_i(n)=W_i((n+1)h)-W_i(nh), \, n=0,...,N-1, \, i=1,2,$$
and $x_i(n)=x_i(nh ,omega).$ $G_i(n)$ are generated using Box-Muller method.

It is shown that the second Euler scheme has the order for weak convergence
1, for sufficiently regular drift and diffusion coefficients.

We assume that $f_i$ in (\ref{aaa1}) are sufficiently smooth such
that the following schemes are well defined.

The second order Euler scheme is defined by the relations
\begin{eqnarray*}
x_i(n+1)&=&x_i(n)+f_i(x(n))h+g_i(x(n))G_i(n)+g_i(x(n))\frac{\partial}{\partial
x_i(n)}g_i(x(n))\frac{G_i(n)^2-h}{2}+\\&+&
\Big[f_i(x(n))\frac{\partial f_i(x(n))}{\partial
x_i(n)}+\frac{1}{2}(g_i(x(n))^2\frac{\partial^2
f_i(x(n))}{\partial x_i(n)\partial
x_i(n)}\Big]\frac{h^2}{2}+\Big[g_i(x(n))\frac{\partial
f_i(x(n))}{\partial
x_i(n)}\\
&+&f_i(x(n))\frac{\partial g_i(x(n))}{\partial
x_i(n)}+\frac{1}{2}(g_i(x(n))^2\frac{\partial^2
g_i(x(n))}{\partial x_i(n)\partial x_i(n)}\Big] \frac{h
G_i(n)}{2},\, i=1,2,
\end{eqnarray*}
where we used the random variables $G_i(n), \, i=1,2.$  In
\cite{Maho}, it is shown that these schemes converge weakly with
order 2.

\end{document}